\documentclass[aop,seceqn,MSNbibl,citesort,dvips]{arximspdf}

\doi{10.1214/10-AOP547}
\volume{39}
\issue{1}
\pubyear{2011}
\firstpage{291}
\lastpage{326}

\makeatletter

\newtheorem{theorem}{Theorem}[section]
\newtheorem{corollary}[theorem]{Corollary}
\newproclaim{definition}[theorem]{Definition}
\newtheorem{lemma}{Lemma}[section]
\newtheorem{prop}[lemma]{Proposition}
\newtheorem{proposition}[theorem]{Proposition}
\newproclaim{remark}[theorem]{Remark}

\newcommand{\RR}{\mathbb{R}}

\newcommand{\ga}{\gamma}

\makeatother

\begin{document}
\begin{frontmatter}

\title{Feynman--Kac formula for heat equation driven by fractional white noise}
\runtitle{Feynman--Kac formula}

\begin{aug}
\author[A]{\fnms{Yaozhong} \snm{Hu}\corref{}\ead[label=e1]{hu@math.ku.edu}},
\author[A]{\fnms{David} \snm{Nualart}\ead[label=e2]{nualart@math.ku.edu}} and
\author[A]{\fnms{Jian} \snm{Song}\ead[label=e3]{jsong@math.ku.edu}}
\runauthor{Y. Hu, D. Nualart and J. Song}
\affiliation{University of Kansas}
\address[A]{Department of Mathematics \\
University of Kansas \\
Lawrence, Kansas, 66045\\
USA\\
\printead{e1}\\
\phantom{E-mail: }\printead*{e2}\\
\phantom{E-mail: }\printead*{e3}} %adresu isvedimo komanda gale!
\end{aug}

\thankstext{t1}{Supported in part by NSF Grant DMS-09-04538.}

% HISTORY:
\received{\smonth{6} \syear{2009}}
\revised{\smonth{2} \syear{2010}}

% ABSTRACT
%
\begin{abstract}
We establish a version of the Feynman--Kac formula for the
multidimensional stochastic heat equation with a multiplicative
fractional Brownian sheet. We use the techniques of Malliavin calculus
to prove that the process defined by the Feynman--Kac formula is a weak
solution of the stochastic heat equation. From the Feynman--Kac
formula, we establish the smoothness of the density of the solution and
the H\"older regularity in the space and time variables. We also derive
a Feynman--Kac formula for the stochastic heat equation in the
Skorokhod sense and we obtain the Wiener chaos expansion of the
solution.
\end{abstract}

% KEYWORDS
%

\pdfkeywords{60H07, 60H15, 60G17, 60G22, 60G30, 35K20, 35R60,
Fractional noise, stochastic heat equations,
Feynman--Kac formula, exponential integrability,
absolute continuity, Holder continuity, chaos expansion}

\end{frontmatter}

%s1 ###
\section{Introduction}

Consider the following heat equation on $\mathbb{R}^d$:
%
%e1.1 ###
%
\begin{equation}\label{e.1.0}
\cases{\dfrac{\partial u}{\partial t}= \dfrac12 \Delta u+c(t,x)u, \vspace*{2pt}\cr
u(0,x)=f(x) ,}
\end{equation}
where $f$ is a bounded measurable function. If $c(t,x)$ is a continuous
function of $(t,x)\in[0,\infty) \times\mathbb{R}^d$, then we have the
well-known Feynman--Kac formula (see \cite{fred}) for the solution of above
equation
\[
u(t,x)=E \biggl[ f(B_{t}^{x})\exp\biggl( \int_{0}^{t}c(t-s,B_{s}^{x})\,ds \biggr) %
\biggr] ,
\]
where $B_{t}^{x}=B_{t}+x$ is a $d$-dimensional Brownian motion starting from
the point $x$.

In this paper, we shall extend the above Feynman--Kac formula to the heat
equation with fractional noise
%
%e1.2 ###
%
\begin{equation}\label{e.1.1}
\cases{\dfrac{\partial u}{\partial t}=\dfrac12 \Delta u+u\,\dfrac
{\partial^{d+1}W}{%
\partial t\,\partial x_1 \cdots\partial x_d}, \vspace*{2pt}\cr
u(0,x)=f(x),}
\end{equation}
where $W(t,x)$ is a fractional Brownian sheet with Hurst parameters $H_{0}$
in time and $(H_1, \ldots, H_d)$ in space, respectively. The difference
between (\ref{e.1.0}) and (\ref{e.1.1}) is that $\frac{\partial
^{d+1}W}{%
\partial t \,\partial x_1 \cdots\partial x_d}$ is no\vspace*{1pt} longer a function
of $t$
and $x$, but a generalized (random) function. For this equation, we can still
formally write down the Feynman--Kac formula
%
%e1.3 ###
%
\begin{equation} \label{feynman}
u(t,x)=E^{B} \biggl[ f(B_{t}^{x})\exp\biggl( \int_{0}^{t}\int_{\mathbb{R}^d}
\delta(B_{t-r}^{x}-y)W(dr,dy) \biggr) \biggr] ,
\end{equation}
where $E^{B}$ denotes the expectation with respect to the Brownian
motion $%
B_{t}^{x}$ and $\delta$ denotes the Dirac delta function.

The aim of this paper is to justify the above formula (\ref{feynman}), to
show that the process $u(t,x)$ is a weak solution of (\ref{e.1.1})
and to establish some properties of this process. First,
we shall show that the stochastic Feynman--Kac functional $%
V_{t,x}:=\int_{0}^{t}\int_{\mathbb{R}^{d}}\delta
(B_{t-r}^{x}-y)W(dr,dy)$ is
a well-defined random variable. This will be done in Section \ref{sec2} using a
suitable approximation of the Dirac delta function, assuming that the Hurst
parameters satisfy $2H_{0}+\sum_{i=1}^{d}H_{i}>d+1$, $H_{0}\geq\frac{1}{2}$
and $H_{i}>\frac{1}{2}$ for $1\leq i\leq d$.

After the definition of the random variable $V_{t,x}$, the next problem is
to show its exponential integrability. With the use of the covariance
structure of the fractional Brownian sheet $W(t,x)$, we show that $u(t,x)$
has exponential moments provided that
%
%e1.4 ###
%
\begin{equation}\label{e.1.4}
E\exp\biggl[ \lambda
\int_{0}^{1}\int_{0}^{1}|r-s|^{2H_{0}-2}%
\prod_{i=1}^{d}|B_{r}^{i}-B_{s}^{i}|^{2H_{i}-2}\,dr\,ds \biggr] <\infty
\end{equation}
for any $\lambda\in\mathbb{R}$. To show that (\ref{e.1.4}) is true,
we use a method introduced by Le Gall in \cite{legall} to derive the
exponential integrability of the renormalized self-intersection local time
of the planar Brownian motion, together with the self-similarity of the
fractional Brownian sheet and several other techniques. This is done in
Section \ref{sec3}.

Another major aim of this paper is to show that $u(t,x)$ defined by %
(\ref{feynman}) is a weak solution of (\ref{e.1.1}). Instead
of following the classical approach based on It\^{o}'s formula, which seems
complicated in our situation, we again use the approximation technique,
together with Malliavin calculus. The main ingredient is to express the
Stratonovich integral as the sum of a Skorokhod integral plus a
correction term
involving Malliavin derivatives. This is a new methodology which is
developed in Section \ref{sec4}.

The Feynman--Kac formula gives an explicit form of a weak solution of
equation (%
\ref{e.1.1}) which turns out to be very useful for obtaining regularity
properties.
Several consequences of this expression are derived in Section \ref{sec5}.
First, we derive the H\"{o}lder continuity of the solution $%
u(t,x)$ with respect to $t$ and $x$, and, afterward, we establish the
smoothness of the density of the probability law of $u(t,x)$ (with respect
to the Lebesgue measure) using techniques of Malliavin calculus.

In the above (\ref{e.1.1}), the solution and the noise are
multiplied using the ordinary product. This gives rise to the Stratonovich
integral when we interpret the equation in its integral form. There are
several papers where the Wick product between the solution and the
noise is
used, this corresponding to the Skorokhod integral. The Stratonovich integral
is more difficult to handle, but it is the right choice if we want to
represent a physical model. Applying a Wiener chaos technique pioneered by
Dawson and Salehi in \cite{DS} and used in several other papers
(see, e.g., the work \cite{hunualart} on the relation between moments of
the solution and self-intersection local times), one can show that there
exists a unique mild solution to the Skorokhod-type equation. We
discuss this
result in Section \ref{sec7} and, using Wiener chaos expansions, we obtain a
Feynman--Kac formula for this solution.

The above techniques work for $H_{i}>1/2$, $i=1,2,\ldots,d$. From the
condition $2H_{0}+\sum_{i=1}^{d}H_{i}>d+1$, it follows that $H_{0}$
must be
greater than $1/2$ and we cannot allow more than one of the $H_{1}
,\ldots
,H_{d}$ to be less than or equal to $1/2$. Thus, if we want to remove the
condition $H_{i}>1/2$, $i=1,2,\ldots,d$, we need $d=1$. We show in
Section \ref{sec7}
that if $d=1$, $H_{1}=\frac{1}{2}$ and $H_{0}>\frac{3}{4}$, then all previous
results hold. When $d=1$, we can also handle the case $H_{0}<1/2$, assuming
that the process has a regular spatial covariance. This has been done
in the
companion paper \cite{HLN}, using different techniques. Finally, the
\hyperref[app]{Appendix} contains some technical results which are used in the
paper.

We would like to close this introduction with some remarks about the
motivation of our work and its connection with other related results.
The existence of a Feynman--Kac
formula like the one we have derived here was mentioned as a conjecture in
a paper by Mocioalca and Viens (see \cite{movi}), although this problem
was circulating long before that. In the lectures by Walsh in Saint
Flour (see \cite{Wa}) it was stated that the one-dimensional equation
in the
It\^o sense driven by a space--time white noise cannot have a
Feynman--Kac formula because the It\^o--Stratonovich correction term is
infinite. In a previous work \cite{hunualart},
two of the present authors considered a Skorokhod-type
equation assuming $H_{i}=\frac{%
1}{2}$ for $i=1,\ldots,d$. In this case, there exists a unique mild
solution obtained by means of the Wiener chaos method if $d=1$ or
$d=2$, $H_{0}>\frac12$ and $t$ is small enough, although the
Feynman--Kac formula is not available unless $d=1$ and $H_{0}>\frac
{3}{4}$ (see Section \ref{sec7}).

A process similar to (\ref{feynman}) was studied by Viens and Zhang in
\cite{VZ}, although it does not have a relation with a stochastic heat
equation and, most likely, the asymptotic results obtained in \cite{VZ}
can be extended to the process (\ref{feynman}).

Recently, Hinz obtained in \cite{Hi} a Feynman--Kac formula for the
stochastic heat equation with a
Gaussian multiplicative noise of the form $\frac{\partial W}{\partial t}(t,x)
$, where $W$ is a fractional Brownian sheet with Hurst parameter $H>%
\frac{1}{2}$ in time and $K\in(0,1)$ in space, and he used this
formula to solve a stochastic Burgers
equation by means of the Hopf--Cole transformation. In this paper, the
noise is more regular in space and this allows the techniques of
classical fractional calculus to be used, together with curvilinear integrals.

%s2 ###
\section{Preliminaries}\label{sec2}

Fix a vector of Hurst parameters $H=(H_0, H_1, \ldots, H_{d})$, where
$H_i\in
( \frac12,1 )$. Suppose that $W=\{W(t,x),t\geq0,x\in\mathbb{R}%
^d\}$ is a zero-mean Gaussian random field with the covariance function%
\[
E(W(t,x)W(s,y))=R_{H_{0}}(s,t)\prod_{i=1}^{d }R_{H_{i}}(x_i, y_i),
\]
where, for any $H\in(0,1)$, we denote by $R_{H}(s,t)$ the covariance
function of the fractional Brownian motion with Hurst parameter $H$,
that is,%
\[
R_{H}(s,t)=\tfrac{1}{2}(|t|^{2H}+|s|^{2H}-|t-s|^{2H}).
\]
In other words, $W$ is a fractional Brownian sheet with Hurst
parameters $%
H_{0}$ in the time variable and $H_{i}$ in the space variables,
$i=1,\ldots,d $.

Denote by $\mathcal{E}$ the linear span of the indicator functions of
rectangles of the form $(s,t]\times(x,y]$ in $\mathbb{R}_{+}\times
\mathbb{R%
}^{d}$. Consider, in $\mathcal{E}$, the inner product defined by%
\[
\bigl\langle I_{(0,s]\times(0,x]},I_{(0,t]\times(0,y]}\bigr\rangle_{\mathcal{H}
}=R_{H_{0}}(s,t)\prod_{i=1}^{d}R_{H_{i}}(x_{i},y_{i}).
\]
In the above formula, if $x_{i}<0$, then we assume, by convention, that
$%
I_{(0,x_{i}]}=-I_{(-x_{i},0]}$. We denote by $\mathcal{H}$ the closure
of $%
\mathcal{E}$ with respect to this inner product. The mapping $%
W\dvtx I_{(0,t]\times(0,x]}\rightarrow W(t,x)$ extends to a linear isometry
between $\mathcal{H}$ and the Gaussian space spanned by $W$. We will denote
this isometry by
\[
W(\phi)=\int_{0}^{\infty}\int_{\mathbb{R}^{d}}\phi(t,x)W(dt,dx)
\]
if $\phi\in\mathcal{H}$. Notice that if $\phi$ and $\psi$ are
functions in $\mathcal{E}$, then%
%
%e2.1 ###
%
\begin{eqnarray}\label{r1}
E ( W(\phi)W(\psi) ) &=& \langle\phi,\psi\rangle_{\mathcal{H}}
\nonumber\\
&=&\alpha_{H}
\int_{\mathbb{R}_{+}^{2}\times\mathbb{R}^{2d}}\phi
(s,x)\psi
(t,y)|s-t|^{2H_{0}-2}\\
&&\hspace*{52.1pt}{}\times\prod_{i=1}^{d}|x_{i}-y_{i}|^{2H_{i}-2}\,ds\,dt\,dx\,dy,\nonumber
\end{eqnarray}
where $\alpha_{H}=\prod_{i=0}^{d}H_{i}(2H_{i}-1)$. Furthermore, $%
\mathcal{H}$ contains the class of measurable functions $\phi$ on $%
\mathbb{R}_{+}\times\mathbb{R}^{d}$ such that%
%
%e2.2 ###
%
\begin{equation} \label{r2}
\int_{\mathbb{R}_{+}^{2}\times\mathbb{R}^{2d}}|\phi(s,x)\phi
(t,y)||s-t|^{2H_{0}-2}\prod
_{i=1}^{d}|x_{i}-y_{i}|^{2H_{i}-2}\,ds\,dt\,dx\,dy<\infty.\hspace*{-32pt}
\end{equation}

We will denote by $D$ the derivative operator in the sense of Malliavin
calculus. That is, if $F$ is a smooth and cylindrical random variable
of the
form%
\[
F=f(W(\phi_{1}),\ldots,W(\phi_{n})),
\]
$\phi_{i}\in\mathcal{H}$, $f\in C_{p}^{\infty}(\mathbb{R}^{n})$
($f$ and
all its partial derivatives have polynomial growth), then $DF$ is the $%
\mathcal{H}$-valued random variable defined by
\[
DF=\sum_{j=1}^{n}\frac{\partial f}{\partial x_{j}}(W(\phi
_{1}),\ldots
,W(\phi_{n}))\phi_{j}.
\]
The operator $D$ is closable from $L^{2}(\Omega)$ into $L^{2}(\Omega;%
\mathcal{H})$ and we define the Sobolev space $\mathbb{D}^{1,2}$ as the
closure of the space of smooth and cylindrical random variables under the
norm%
\[
\| DF \| _{1,2}=\sqrt{E(F^{2})+E( \| DF \| _{\mathcal{H}%
}^{2})}.
\]
We denote by $\delta$ the adjoint of the derivative operator given by
duality formula
%
%e2.3 ###
%
\begin{equation}\label{dua}
E(\delta(u)F)=E ( \langle DF,u \rangle_{\mathcal{H}} )
\end{equation}
for any $F\in\mathbb{D}^{1,2}$ and any element $u\in$ $L^{2}(\Omega;
\mathcal{H})$ in the domain of $\delta$. The operator $\delta$ is also
called the \textit{Skorokhod integral} because in the case of the
Brownian motion, it
coincides with an extension of the It\^{o} integral introduced by Skorokhod.
We refer to Nualart \cite{nualart} for a detailed account of the Malliavin
calculus with respect to a Gaussian process.
If $DF$ and $u$ are almost surely measurable functions on $\mathbb{R}_+
\times\mathbb{R}^d$ verifying
condition (\ref{r2}), then the duality formula (\ref{dua}) can be
written using the expression of the inner product in $\mathcal{H}$
given in (\ref{r1}):
\begin{eqnarray*}
E ( \delta(u)F ) &=& \alpha_H
E \biggl( \int_ {\mathbb{R}^2_+\times\mathbb{R}^{2d}}
D_{s,x} F u(t,y) |s-t|^{2H_{0}-2}\\
&&\hspace*{66.2pt}{}\times\prod
_{i=1}^{d}|x_{i}-y_{i}|^{2H_{i}-2}\,ds\,dt\,dx\,dy \biggr).
\end{eqnarray*}
We recall the following
formula, which we will use in the paper:%
%
%e2.4 ###
%
\begin{equation}\label{a5}
FW(\phi)=\delta(F\phi)+ \langle DF,\phi\rangle_{\mathcal{H}}
\end{equation}
for any $\phi\in\mathcal{H}$ and any random variable $F$ in the Sobolev
space $\mathbb{D}^{1,2}$.

Throughout the paper, $C$ will denote a positive constant which may
vary from one
formula to another.

%s3 ###
\section{Definition and exponential integrability of the stochastic
Feynman--Kac functional}\label{sec3}

For any $\varepsilon>0$, we denote by $p_{\varepsilon}(x)$ the
$d$-dimensional heat kernel:
\[
p_{\epsilon}(x)=(2\pi\varepsilon)^{-{d}/{2}}e^{-{|x|^{2}}/{
2\epsilon}},\qquad x\in\mathbb{R}^{d}.
\]
On the other hand, for any $\delta>0$, we define the function
\[
\varphi_{\delta}(x)=\frac{1}{\delta}I_{[0,\delta]}(x).
\]
$\varphi_{\delta}(t)p_{\varepsilon}(x)$ then provides an approximation
of the Dirac delta function $\delta(t,x)$ as $\varepsilon$ and
$\delta$
tend to zero. We denote by $W^{\epsilon,\delta}$ the approximation of the
fractional Brownian sheet $W(t,x)$ defined by
%
%e3.1 ###
%
\begin{equation} \label{e1}
W^{\epsilon,\delta}(t,x)=\int_{0}^{t}\int_{\mathbb{R}^{d}}\varphi
_{\delta
}(t-s)p_{\epsilon}(x-y)W(s,y)\,ds\,dy .
\end{equation}

Fix $x\in\mathbb{R}^{d}$ and $t>0$. Suppose that $B=\{B_{t},t\geq0\}$
is a
$d$-dimensional standard Brownian motion independent of $W$. We denote
by $%
B_{t}^{x}=B_{t}+x$ the Brownian motion starting at the point $x$. We are
going to define the random variable $\int_{0}^{t}\int_{\mathbb
{R}^{d}}\delta
(B_{t-r}^{x}-y)W(dr,dy)$ by approximating the Dirac delta function
$\delta
(B_{t-r}^{x}-y)$ by%
%
%e3.2 ###
%
\begin{equation} \label{e2}
A_{t,x}^{\epsilon,\delta}(r,y)=\int_{0}^{t}\varphi_{\delta
}(t-s-r)p_{\varepsilon}(B_{s}^{x}-y)\,ds.
\end{equation}
We will show that for any $\varepsilon>0$ and $\delta>0$, the function
$A_{t,x}^{\epsilon,\delta}$ belongs to the space $\mathcal{H}$ almost
surely and the family of random variables%
%
%e3.3 ###
%
\begin{equation} \label{e3}
V_{t,x}^{\varepsilon,\delta}= \int_{0}^{t}\int_{\mathbb{R}%
^{d}}A_{t,x}^{\epsilon,\delta}(r,y)W(dr,dy)
\end{equation}
converges in $L^{2}$ as $\varepsilon$ and $\delta$ tend to zero.

The specific approximation chosen here will allow us, in Section \ref{sec4}, to
construct an approximate
Feynman--Kac formula with the random potential $\dot{W}^{\epsilon
,\delta}(t,x)$ given in
(\ref{pot}). Moreover, this approximation has the useful properties
proved in
Lem\-mas~\ref{lemma2} and \ref{lemma3}. We could have used other types of
approximation schemes with similar results. Also, we can restrict
ourselves to the special case ${\delta}={\varepsilon}$, but the slightly
more general case considered here does not need any additional effort.

Throughout the paper, we denote by $E^B(\Phi(B,W))$ [resp., by
$E^W(\Phi(B,W))$]
the expectation of a functional $\Phi(B,W)$ with respect to $B$
(resp., with
respect to~$W$). We will use $E$ for the composition $E^BE^W$ and also
in the
case of a random variable depending only on $B$ or $W$.
\begin{theorem}
\label{teo1}
Suppose that $2H_{0}+\sum_{i=1}^d H_{i}>d+1$. Then, for any
$%
\varepsilon>0 $ and $\delta>0$, $A_{t,x}^{\varepsilon,\delta}$
defined in
(\ref{e2}) belongs to $\mathcal{H}$ and the family of random
variables $
V_{t,x}^{\epsilon,\delta}$ defined in (\ref{e3}) converges in $L^{2}$
to a
limit denoted by (this being the stochastic Feynman--Kac functional)
%
%e3.4 ###
%
\begin{equation}\label{e.3.4}
V_{t,x}=\int_{0}^{t}\int_{\mathbb{R}^d}\delta
(B_{t-r}^{x}-y)W(dr,dy) .
\end{equation}
Conditional on $B$, $V_{t,x}$ is a Gaussian random variable with mean $0$
and variance%
%
%e3.5 ###
%
\begin{equation}\label{e.2.14}
\operatorname{Var}^W(V_{t,x}) =\alpha_H \int_{0}^{t}\int_{0}^{t}|r -s|^{2H_{0}-2}
\prod_{i=1}^d | B^i_{r }-B^i_{s} | ^{2H_{i}-2}\,dr \,ds .
\end{equation}
\end{theorem}
\begin{pf}
Fix $\varepsilon$, $\varepsilon^{\prime}$, $\delta$ and $\delta
^{\prime
}>0$. Let us compute the inner product%
%
%e3.6 ###
%
\begin{eqnarray}\label{e4}
&&
\langle A_{t,x}^{\varepsilon,\delta},A_{t,x}^{\varepsilon^{\prime
},\delta^{\prime}} \rangle_{\mathcal{H}}\nonumber\\
&&\qquad=\alpha
_{H}\int_{[0,t]^{4}}\int_{\mathbb{R}^{2d}}p_{\epsilon
}(B_{s}^{x}-y)p_{\epsilon^{\prime}}(B_{r}^{x}-z)\nonumber\\[-8pt]\\[-8pt]
&&\qquad\quad\hspace*{60pt}{}\times\varphi_{\delta}(t-s-u)\varphi_{\delta^{\prime}}(t-r-v)
|u-v|^{2H_{0}-2}\nonumber\\
&&\qquad\quad\hspace*{60pt}{}\times\prod_{i=1}^{d}|y_{i}-z_{i}|^{2H_{i}-2}\,dy\,dz\,du\,dv\,ds\,dr.\nonumber
\end{eqnarray}
By Lemmas \ref{lemma2} and \ref{lemma3}, we have the estimate
%
%e3.7 ###
%
\begin{eqnarray} \label{e5}
&&\int_{[0,t]^{2}}\int_{\mathbb{R}^{2d}}p_{\epsilon
}(B_{s}^{x}-y)p_{\epsilon^{\prime}}(B_{r}^{x}-z) \varphi_{\delta}(t-s-u)\varphi_{\delta^{\prime
}}(t-r-v)\nonumber\\
&&\hspace*{34pt}\quad{}\times|u-v|^{2H_{0}-2}\prod_{i=1}^{d}|y_{i}-z_{i}|^{2H_{i}-2}\,dy\,dz\,du\,dv
\\
&&\hspace*{36.5pt}\qquad\leq C|s-r|^{2H_{0}-2}\prod_{i=1}^{d} | B_{s}^{i}-B_{r}^{i} |
^{2H_{i}-2}\nonumber
\end{eqnarray}
for some constant $C>0$. The expectation of this random variable is
integrable in $[0,t]^{2}$ because%
%
%e3.8 ###
%
\begin{eqnarray} \label{3.8}
&&E^{B}\int_{0}^{t}\int_{0}^{t}|s-r|^{2H_{0}-2}\prod_{i=1}^{d} |
B_{s}^{i}-B_{r}^{i} | ^{2H_{i}-2}\,ds\,dr \nonumber\\
&&\qquad=\prod_{i=1}^{d}E|\xi
|^{2H_{i}-2}\int_{0}^{t}\int_{0}^{t}|s-r|^{2H_{0}+%
\sum_{i=1}^{d}H_{i}-d-2}\,ds\,dr \\
&&\qquad=\frac{2\prod_{i=1}^{d}E|\xi|^{2H_{i}-2}t^{\kappa+1}}{\kappa(
\kappa+1 ) } <\infty,\nonumber
\end{eqnarray}
where
%
%e3.9 ###
%
\begin{equation} \label{kappa}
\kappa=2H_{0}+\sum_{i=1}^{d}H_{i}-d-1>0
\end{equation}
and $\xi$ is a $N(0,1)$ random variable.

As a consequence, taking the mathematical expectation with respect to
$B$ in
(\ref{e4}), letting $\varepsilon=\varepsilon^{\prime}$ and $
\delta=\delta^{\prime}$ and using the estimates (\ref{e5}) and
(\ref
{3.8}%
) yields%
\[
E^{B} \| A_{t,x}^{\varepsilon,\delta} \| _{\mathcal{H}}^{2}\leq C.
\]
This implies that almost surely $A_{t,x}^{\varepsilon,\delta}$
belongs to
the space $\mathcal{H}$ for all $\varepsilon$ and $\delta>0$. Therefore,
the random variables $V_{t,x}^{\varepsilon,\delta
}=W(A_{t,x}^{\varepsilon
,\delta})$ are well defined and we have
\[
E^{B}E^{W}(V_{t,x}^{\epsilon,\delta}V_{t,x}^{\epsilon^{\prime
},\delta
^{\prime}})=E^{B} \langle A_{t,x}^{\varepsilon,\delta
},A_{t,x}^{\varepsilon^{\prime},\delta^{\prime}} \rangle_{\mathcal{
H}}.
\]
For any $s\neq r$ and $B_{s}\neq B_{r}$, as $\varepsilon$,
$\varepsilon
^{\prime}$, $\delta$ and $\delta^{\prime}$ tend to zero, the left-hand
side of the inequality (\ref{e5}) converges to $|s-r|^{2H_{0}-2}%
\prod_{i=1}^{d} | B_{s}^{i}-B_{r}^{i} | ^{2H_{i}-2}$. Therefore, by
the dominated convergence theorem, we obtain that
$E^{B}E^{W}(V_{t,x}^{\epsilon
,\delta}V_{t,x}^{\epsilon^{\prime},\delta^{\prime}})$ converges to
$%
\Sigma_{t}$ as $\varepsilon$, $\varepsilon^{\prime}$, $\delta$ and
$%
\delta^{\prime}$ tend to zero, where
\[
\Sigma_{t}=\frac{2\alpha_{H}\prod_{i=1}^{d}E|\xi
|^{2H_{i}-2}t^{\kappa+1}%
}{\kappa( \kappa+1 ) }.
\]
Thus, we obtain
\[
E ( V_{t,x}^{{\varepsilon},{\delta}}-V_{t,x}^{{\varepsilon}^{\prime
},{\delta}^{\prime} } ) ^{2}=E ( V_{t,x}^{{\varepsilon},{%
\delta} } ) ^{2}-2E ( V_{t,x}^{{\varepsilon},{\delta}%
}V_{t,x}^{{\varepsilon}^{\prime},{\delta}^{\prime}} ) +E (
V_{t,x}^{{\varepsilon}^{\prime},{\delta}^{\prime} } )
^{2}\rightarrow0 .
\]
This implies that $V_{t,x}^{\epsilon_{n},\delta_{n}}$ is a Cauchy
sequence in $L^{2}$ for all sequences $\varepsilon_{n}$ and $\delta_{n}$
converging to zero. As a consequence, $V_{t,x}^{\epsilon_{n},\delta_{n}}$
converges in $L^{2}$ to a limit denoted by $V_{t,x}$ which does not depend
on the choice of the sequences $\varepsilon_{n}$ and $\delta_{n}$.
Finally, by a similar argument, we show (\ref{e.2.14}).
\end{pf}

Condition $2H_{0}+\sum_{i=1}^{d}H_{i}> d+1$ is sharp and cannot be
improved. In fact, if this condition does not hold, then almost surely $(r,y)
\mapsto\delta(B^x_{t-r} -y)$ is not an element of the space $\mathcal
{H}$, as follows from the next proposition.
\begin{proposition}
Suppose that $H_i>1/2$, $i=0, 1, \ldots, d$, and $2H_{0}+\sum
_{i=1}^{d}H_{i}\le d+1$.
Then, conditionally on $B$, the family $V_{t,x}^{{\varepsilon},{\delta
}}$ does not converge in probability as $\varepsilon$ and $\delta$ tend
to zero for almost all trajectories of $B$.
\end{proposition}
\begin{pf}
Given $B$, $V_{t,x}^{{\varepsilon},{\delta}}$ is a Gaussian family of
random variables and it suffices to show that they do not converge in
$L^2$. This follows from the fact that the variance limit is infinite
almost surely. In fact, from the L\'evy modulus of continuity of the
Brownian motion, it is easy to show that if $2H_{0}+\sum
_{i=1}^{d}H_{i}\le d+1$, then
\[
\int_0^t\int_0^t |s-r|^{2H_{0}-2}\prod_{i=1}^{d} |
B_{s}^{i}-B_{r}^{i} |
^{2H_{i}-2} \,ds\,dr =\infty
\]
almost surely.
\end{pf}

The next result provides the exponential integrability of the random
variable $V_{t,x}$ defined in (\ref{e.3.4}).
\begin{theorem}
\label{teo2} Suppose that $2H_0+ \sum_{i=1}^d H_{i}>d+1$. Then, for any
$%
\lambda\in\mathbb{R}$, we have
%
%e3.10 ###
%
\begin{equation} \label{e21}
E \exp\biggl( \lambda\int_{0}^{t}\int_{\mathbb{R}^d}\delta
(B_{t-r}^{x}-y)W(dr,dy) \biggr) <\infty.
\end{equation}
\end{theorem}
\begin{pf}
The proof involves several steps.

\textit{Step} 1. From (\ref{e.2.14}), we obtain
\[
Ee^{\lambda V_{t,x}}=E^{B}\exp\biggl( \frac{\lambda^{2}}{2}\alpha
_{H}\int_{0}^{t}\int_{0}^{t}|s-r|^{2H_{0}-2}%
\prod_{i=1}^{d}|B_{s}^{i}-B_{r}^{i}|^{2H_{i}-2}\,ds\,dr \biggr)
\]
and the scaling property of the Brownian motion yields
%
%e3.11 ###
%
\begin{equation} \label{e9a}
Ee^{\lambda V_{t,x}}=Ee^{\mu Y},
\end{equation}
where $\mu=\frac{\lambda^{2}}{2}\alpha_{H}t^{\kappa+1}$, $\kappa$
is as defined in (\ref{kappa}) and
%
%e3.12 ###
%
\begin{equation} \label{e9}
Y=\int_{0}^{1}\int_{0}^{1}|s-r|^{2H_{0}-2}%
\prod_{i=1}^{d}|B_{s}^{i}-B_{r}^{i}|^{2H_{i}-2}\,ds\,dr.
\end{equation}
It then suffices to show that the random variable $Y$ has exponential
moments of all orders.

\textit{Step} 2. Our approach to proving that $E\exp(
\lambda Y) <\infty$ for any $\lambda\in\mathbb{R}$ is motivated by the
method of Le Gall \cite{legall}. For $k=1,\ldots,2^{n-1}$, we define $
A_{n,k}= [ \frac{2k-2}{2^{n}},\frac{2k-1}{2^{n}} ] \times[
\frac{2k-1}{2^{n}},\frac{2k}{2^{n}} ] $ and
\[
\alpha_{n,k}=\int_{A_{n,k}}|s-r|^{2H_{0}-2} \prod_{i=1}^d
|B^i_{s}-B^i_{r}|^{2H_i-2}\,ds\,dr.
\]
The random variables $\alpha_{n,k}$ have the following two properties:

\begin{longlist}
\item for every $n\ge1$, the variables $\alpha_{n,1},\ldots,%
\alpha_{n,2^{n-1}}$ are independent;

\item $\alpha_{n,k}\stackrel{d}{=}2^{-n(\kappa+1)} \alpha
_0$, where
\[
\alpha_0= \int_{0}^{1}\int_{0}^{1}(s+r)^{2H_{0}-2} \prod_{i=1}^d
|B^i_{s}-\tilde B^i_{r}|^{2H_i-2}\,ds\,dr,
\]
and $\tilde B $ is a standard Brownian motion independent of $B$.
\end{longlist}

The condition $2H_{0}+\sum_{i=1}^{d}H_{i}>d+1$ implies that $E\alpha
_{0}<\infty$ and we deduce that
\[
Y=2\sum_{n=1}^{\infty}\sum_{k=1}^{2^{n-1}}\alpha_{n,k},
\]
where the series converges in the $L^{1}$ sense.

\textit{Step} 3. For any integer $n\geq1$, we
claim that
%
%e3.13 ###
%
\begin{equation} \label{e31}
E\alpha_{0}^{n}\leq E \Biggl( C\int_{0}^{1} \prod_{i=1}^d |B^i_{s}
|^{2H_i-2}\,ds \Biggr) ^{n}
\end{equation}
for some constant $C>0$. In fact, we have%
%
%e3.14 ###
%
\begin{equation}\label{e41}
E \alpha_{0} ^{n}=E\int_{[0,1]^{2n}}\prod_{j=1}^{n} \prod_{i=1}^d (
|s_{j}+t_{j}|^{2H_{0}-2} |B^i_{s_{j}}-\tilde B^i_{t_{j}} | ^{ 2H_i-2})
\,ds\,dt.
\end{equation}
Using the formula
\[
c^{-z}=\frac{1}{\Gamma(z)}\int_{0}^{\infty}e^{-c\tau}\tau
^{z-1}\,d\tau,
\]
we obtain, for each $i=1,\ldots,d$,
%
%e3.15 ###
%
\begin{eqnarray} \label{e6}
E\prod_{j=1}^{n} |B^i_{s_{j}}-\tilde B^i_{t_{j}} | ^{2H_i-2}& =&
\Gamma(
1-H_i)^{-n} \nonumber\\
&& \times\int_{[0,\infty)^{n}}E\exp\Biggl(
-\sum_{j=1}^{n}|B^i_{s_{j}}-\tilde B^i_{t_{j}} |^{2}\tau_{j} \Biggr)\\
&&\hspace*{39.7pt}{}\times\prod_{j=1}^{n}\tau_{j}^{-H_i}\,d\tau.\nonumber
\end{eqnarray}
For any $\tau_{1},\ldots,\tau_{n}>0$ and $s_{1},t_{1},\ldots
,s_{n},t_{n}\in(0,1)$, we define
\[
Q_{1}= \bigl( E(B^i_{s_{j}}B^i_{s_{k}})\sqrt{\tau_{j}\tau_{k}} \bigr)
_{n\times n},\qquad Q_{2}= \bigl( E(\tilde B^i_{t_{j}} \tilde B^i_{t_{k}} )%
\sqrt{\tau_{j}\tau_{k}} \bigr) _{n\times n}.
\]
We know that%
%
%e3.16 ###
%
\begin{equation} \label{e51}
E\exp\Biggl( -\sum_{j=1}^{n}|B^i_{s_{j}}-\tilde B^i_{t_{j}} |^{2}\tau
_{j} \Biggr) =\det(I+2Q_{1}+2Q_{2})^{-{1/2}} .
\end{equation}
Substituting (\ref{e51}) into (\ref{e6}) yields
%
%e3.17 ###
%
\begin{eqnarray} \label{e7}\quad
&&E\prod_{j=1}^{n}|B^i_{s_{j}}-\tilde B^i_{t_{j}} |^{2H_i-2}\nonumber\\
&&\qquad=\Gamma
(1-H_i)^{-n}
\int_{[0,\infty)^{n}}\det(I+2Q_{1}+2Q_{2})^{-{1/2}%
}\prod_{j=1}^{n}\tau_{j}^{-H_i}\,d\tau\nonumber\\
&&\qquad\leq\Gamma(1-H_i)^{-n}\nonumber\\
&&\qquad\quad{}\times\int_{[0,\infty)^{n}}\det
(I+2Q_{1})^{-
{1/4}%
}\det(I+2Q_{2})^{-{1/4}}\prod_{j=1}^{n}\tau_{j}^{ -H_i}\,d\tau
\\
&&\qquad\leq\Gamma(1-H_i)^{-n} \Biggl[ \int_{[0,\infty)^{n}}\det(I+2Q_{1})^{-%
{1/2}}\prod_{j=1}^{n}\tau_{j}^{-H_i}d\tau\Biggr] ^{{1}/{2}}
\nonumber\\
&&\qquad\quad{}\times\Biggl[ \int_{[0,\infty)^{n}}\det(I+2Q_{2})^{-{1/2}%
}\prod_{j=1}^{n}\tau_{j}^{-H_i}\,d\tau\Biggr] ^{{1}/{2}} \nonumber\\
&&\qquad= \Biggl[ E\prod_{j=1}^{n}|B^i_{s_{j}}|^{2H_i-2}E\prod_{j=1}^{n}|\tilde
B^i_{t_{j}} |^{2H_i-2} \Biggr] ^{{1}/{2}} ,\nonumber
\end{eqnarray}
where, in the above first inequality, we have used the estimates
\begin{eqnarray*}
(I+2Q_{1}+2Q_{2})&\geq&\tfrac{1}{2}[(I+2Q_{1})+(I+2Q_{2})]\\
&\geq&
(I+2Q_{1})^{{1}/{2}}(I+2Q_{2})^{{1}/{2}}.
\end{eqnarray*}
Substituting (\ref{e7}) into (\ref{e41}) and using the inequality $%
(s_j+t_j)^ {2H_{0}-2} \le s_j^{H_0-1}\times t_j^{H_0-1}$, we obtain
\begin{eqnarray*}
E\alpha_{0}^{n} &\leq
&\int_{[0,1]^{2n}}\prod_{j=1}^{n}(s_{j}+t_{j})^{2H_{0}-2} \prod_{i=1}^d
\Biggl[ E\prod_{j=1}^{n}|B^i_{s_{j}}|^{2H_i-2 }E\prod_{j=1}^{n}|\tilde
B^i_{t_{j}} |^{2H_i-2} \Biggr] ^{{1}/{2}}\,ds\,dt \\
&\le& \Biggl( \int_{[0,1]^{n}} \prod_{j=1}^{n} s_j^{ H_0-1} \Biggl(
E\prod_{j=1}^{n} \prod_{i=1}^d |B^i_{s_{j}}|^{2H_i-2} \Biggr) ^{{1}/{2
}%
}\,ds \Biggr) ^{2 } .
\end{eqnarray*}
Finally, using H\"older's inequality with $\frac1 {H_0} <p <2$, we get
\begin{eqnarray*}%\label{e.3.15}
E\alpha_{0}^{n} &\leq& C^n \Biggl( \int_{[0,1]^{n}} \Biggl( E
\prod_{i=1}^d\prod_{j=1}^{n} |B^i_{s_{j}}|^{2H_i-2 } \Biggr)^{p/2} \,ds
\Biggr) ^{2/p} \\
&\leq&C^{n} \int_{[0,1]^{n}} E \prod_{i=1}^d\prod_{j=1}^{n}
|B^i_{s_{j}}|^{2H_i-2} \,ds \\
&=& E \Biggl( C\int_{0}^{1} \prod_{i=1}^d |B^i_{s}|^{2H_i-2}\,ds \Biggr) ^{n}.
\end{eqnarray*}
This completes the proof of (\ref{e31}).

\textit{Step} 4. For any $\lambda>0$, using (\ref{e31})
and Lemma \ref{lemma5} in the \hyperref[app]{Appendix}, we obtain
%
%e3.18 ###
%
\begin{equation}\label{e8}
Ee^{\lambda\alpha_{0}}\leq E\exp\Biggl( C \lambda\int_{0}^{1}
\prod_{i=1}^d |B^i_{s}|^{2H_i-2}\,ds \Biggr) <\infty,
\end{equation}
because $\rho<1$.

\textit{Step} 5. Define $\varphi(\lambda
)=E(e^{\lambda
(\alpha_{0}-E\alpha_{0})})$. By (\ref{e8}), $\varphi(\lambda
)<\infty$
for all $\lambda\in\mathbb{R}$. Since $\varphi^{\prime}(0)=0$, for every
$K>0$, we can find a positive constant $C_{K}$ such that for all
$\lambda\in
[0,K]$,
\[
\varphi(\lambda)\leq1+C_{K}\lambda^{2}.
\]
Define $\overline{\alpha}_{n,k}={\alpha}_{n,k}-E({\alpha}_{n,k})$.
Fix $%
K>0$ and $a\in(0,\kappa+1)$, where $\kappa$ is as defined in (\ref%
{kappa}). Recall that by property (ii) in step 3, $\overline{\alpha
}_{n,k}%
\stackrel{d}{=}2^{-n(\kappa+1)}\overline{\alpha}_{0}$. For every
$N\geq2$,
set $b_{N}=2K\prod_{j=2}^{j=N}(1-2^{-a(j-1)})$ and $b_{1}=2K$. Then, by
H%
\"{o}lder's inequality and properties (i) and (ii) of ${\alpha
}_{n,k}$, we
have, for $N\geq2$,
\begin{eqnarray*}
&& E\exp\Biggl( b_{N}\sum_{n=1}^{N}\sum_{k=1}^{2^{n-1}}\overline{\alpha}%
_{n,k} \Biggr) \\
&&\qquad\leq \Biggl[ E\exp\Biggl( \frac{b_{N}}{1-2^{-a(N-1)}}\sum_{n=1}^{N-1}%
\sum_{k=1}^{2^{n-1}}\overline{\alpha}_{n,k} \Biggr) \Biggr]
^{1-2^{-a(N-1)}} \\
&&\qquad\quad{} \times\Biggl[ E\exp\Biggl( 2^{a(N-1)}b_{N}\sum_{k=1}^{2^{N-1}}\overline{%
\alpha}_{N,k} \Biggr) \Biggr] ^{2^{-a(N-1)}} \\
&&\qquad \leq E\exp\Biggl( b_{N-1}\sum_{n=1}^{N-1}\sum_{k=1}^{2^{n-1}}\overline{%
\alpha}_{n,k} \Biggr) \varphi\bigl(b_{N}2^{a(N-1)-(\kappa
+1)N}\bigr)^{2^{(1-a)(N-1)}}.
\end{eqnarray*}
Notice that $b_{N}2^{a(N-1)-(\kappa+1)N}\leq2K$. It follows that
\begin{eqnarray*}
\varphi\bigl(b_{N}2^{a(N-1)-(\kappa+1)N}\bigr)^{2^{(1-a)(N-1)}}& \leq &\bigl(
1+C_{K}b_{N}^{2}2^{2((a-\kappa-1)N-a)} \bigr) ^{2^{(1-a)(N-1)}} \\
&\leq&\exp\bigl(C2^{(a+1-2(\kappa+1))N}\bigr)
\end{eqnarray*}
for a constant $C$ independent of $N$. By induction, we get
\begin{eqnarray*}
E\exp\Biggl( b_{N}\sum_{n=1}^{N}\sum_{k=1}^{2^{n-1}}\overline{\alpha}%
_{n,k} \Biggr) & \leq &\exp\Biggl( C\sum_{n=2}^{N}2^{(a+1-2(\kappa
+1))n} \Biggr) E\exp( b_{1}\overline{\alpha}_{1,1} ) \\
& \leq &\exp\bigl( C\bigl(1-2^{a+1-2(\kappa+1)}\bigr)^{-1} \bigr) \varphi(K).
\end{eqnarray*}
Letting $N$ tend to infinity and using Fatou's lemma, we obtain
\[
E\exp\bigl( b_{\infty}(Y-EY)/2 \bigr) <\infty,
\]
where $b_{\infty}=2K\prod_{j=1}^{\infty}(1-2^{-aj})>0$. Since $K>0$ is
arbitrary, we conclude that $E\exp(\lambda Y)<\infty$ for all
$\lambda\in
\mathbb{R}$. This completes the proof, in view of (\ref{e9a}).
\end{pf}

%s4 ###
\section{Feynman--Kac formula}\label{sec4}

We recall that $W$ is a fractional Brownian sheet on $\mathbb
{R}_{+}\times
\mathbb{R}^{d}$ with Hurst parameters $(H_{0},H_{1},\ldots,H_{d})$,
where $%
H_{i}\in(\frac{1}{2},1)$ for $i=0,\ldots,d$. For any $\varepsilon$, $
\delta>0$, we define%
%
%e4.1 ###
%
\begin{equation} \label{pot}
\dot{W}^{\epsilon,\delta}(t,x):=\int_{0}^{t}\int_{\mathbb
{R}^{d}}\varphi
_{\delta}(t-s)p_{\epsilon}(x-y)W(ds,dy).
\end{equation}
In order to provide a notion of solution for the heat equation with fractional
noise~(\ref{e.1.1}), we need the following definition of the Stratonovich
integral, which is equivalent to that of Russo and Vallois in \cite{RV}.
\begin{definition}
\label{def2} Given a random field $v=\{v(t,x),t\geq0,x\in\mathbb
{R}^d\}$
such that
\[
\int_{0}^{T}\int_{\mathbb{R}^d}|v(t,x)|\,dx\,dt<\infty
\]
almost surely for all $T>0$, the Stratonovich integral $\int
_{0}^{T}\int
_{%
\mathbb{R}^d}v(t,x)W(dt,dx)$ is defined as the following limit in
probability, if it exists:%
\[
\lim_{\epsilon,\delta\downarrow0}\int_{0}^{T}\int_{\mathbb
{R}^d}v(t,x)%
\dot{W}^{\epsilon,\delta}(t,x)\,dx\,dt.
\]
\end{definition}

We are going to consider the following notion of solution for
(\ref%
{e.1.1}).
\begin{definition}
\label{def} A random field $u=\{u(t,x),t\geq0,x\in\mathbb{R}^d\}$
is a
weak solution of (\ref{e.1.1}) if, for any $C^{\infty}$
function $%
\varphi$ with compact support on $\mathbb{R}^d$, we have
\begin{eqnarray*}
\int_{\mathbb{R}^d}u(t,x)\varphi(x)\,dx&=&\int_{\mathbb
{R}^d}f(x)\varphi
(x)\,dx + \frac12 \int_{0}^{t}\int_{\mathbb{R}^d}u(s,x) \Delta
\varphi(x)\,dx\,ds
\\
&&{} +\int_{0}^{t}\int_{\mathbb{R}^d}u(s,x)\varphi(x)W(ds,dx)
\end{eqnarray*}
almost surely for all $t\ge0$, where the last term is a Stratonovich
stochastic integral in the sense of Definition \ref{def2}.
\end{definition}

The following is the main result of this section.
\begin{theorem} \label{th1}
Suppose that $2H_0+\sum_{i=1}^d H_{i}>d+1$ and that $f$ is a bounded
measurable function. Then, the process%
%
%e4.2 ###
%
\begin{equation} \label{utx}
u(t,x)=E^{B} \biggl( f(B_{t}^{x})\exp\biggl( \int_{0}^{t}\int_{\mathbb{R}%
^d}\delta(B_{t-r}^{x}-y)W(dr,dy) \biggr) \biggr)
\end{equation}
is a weak solution of (\ref{e.1.1}).
\end{theorem}
\begin{pf}
Consider the approximation of (\ref{e.1.1}) given by the
following heat equation with a random potential:
%
%e4.3 ###
%
\begin{equation}\label{approx}
\cases{\dfrac{\partial u^{\varepsilon,\delta}}{\partial t}= \dfrac12
\Delta
u^{\varepsilon,\delta}+u^{\varepsilon,\delta}\dot
{W}_{t,x}^{\epsilon
,\delta}, \vspace*{2pt}\cr
u^{\varepsilon,\delta}(0,x)=f(x).}
\end{equation}
From the classical Feynman--Kac formula, we know that%
\[
u^{\varepsilon,\delta}(t,x) =E^{B} \biggl( f(B_{t}^{x})\exp\biggl(
\int_{0}^{t}\dot{W}^{\epsilon,\delta}(t-s,B_{s}^{x})\,ds \biggr) \biggr),
\]
where $B^x_t$ is a $d$-dimensional Brownian motion independent of $W$
starting at $x$. By Fubini's theorem,
we can write%
\begin{eqnarray*}
\int_{0}^{t}\dot{W}^{\varepsilon,\delta}(t-s,B_{s}^{x})\,ds
&=&\int_{0}^{t} \biggl( \int_{0}^{t}\int_{\mathbb{R}^d}\varphi_{\delta
}(t-s-r)p_{\varepsilon}(B_{s}^{x}-y)W(dr,dy) \biggr) \,ds \\
&=&\int_{0}^{t} \int_{\mathbb{R}^d} \biggl( \int_{0}^{t}\varphi_{\delta
}(t-s-r)p_{\varepsilon}(B_{s}^{x}-y)\,ds \biggr) W(dr,dy) \\
&=&V_{t,x}^{\varepsilon,\delta},
\end{eqnarray*}
where $V_{t,x}^{\varepsilon,\delta}$ is defined in (\ref{e3}).
Therefore,%
\[
u ^{\varepsilon,\delta}(t,x)=E^{B} ( f(B_{t}^{x})\exp(
V_{t,x}^{\varepsilon,\delta} ) ) .
\]

\textit{Step} 1. We will prove that for any $x\in
\mathbb{R%
}^d$ and any $t>0$, we have%
%
%e4.4 ###
%
\begin{equation} \label{a1}
\lim_{\varepsilon,\delta\downarrow0} E^W |u ^{\varepsilon,\delta
}(t,x)-u(t,x)|^{p}=0
\end{equation}
for all $p\geq2$, where $u(t,x)$ is defined in (\ref{utx}). Notice that
\begin{eqnarray*}
E^W |u ^{\varepsilon,\delta}(t,x)-u(t,x)|^{p} &=&E^{W} \bigl| E^{B} \bigl(
f(B_{t}^{x}) [ \exp( V_{t,x}^{\varepsilon,\delta} ) -\exp
( V_{t,x} ) ] \bigr) \bigr| ^{p} \\
&\leq& \| f \| _{\infty}^{p}E | \exp(
V_{t,x}^{\varepsilon,\delta} ) -\exp( V_{t,x} ) |
^{p},
\end{eqnarray*}
where $V_{t,x}$ is defined in (\ref{e.3.4}). Since $\exp(
V_{t,x}^{\varepsilon,\delta} ) $ converges to $\exp(
V_{t,x} ) $ in probability by Theorem \ref{teo1}, to show (\ref{a1}), it
suffices to prove that for any $\lambda\in\mathbb{R}$,
%
%e4.5 ###
%
\begin{equation} \label{a2}
\sup_{\epsilon,\delta}E\exp( \lambda V_{t,x}^{\varepsilon,\delta
} ) <\infty.
\end{equation}
The estimate (\ref{a2}) follows from (\ref{e3}), (\ref{e5}) and
(\ref
{e21}%
):
%
%e4.6 ###
%
\begin{eqnarray} \label{a2a}\hspace*{8pt}
E \exp( \lambda V_{t,x}^{\varepsilon,\delta} ) &=& E \exp
\biggl( \frac{\lambda^{2}}{2} \| A_{t,x}^{\varepsilon,\delta} \|
_{\mathcal{H}}^{2} \biggr) \nonumber\\
&\leq& E \exp\Biggl( \frac{\lambda^{2}}{2}C\int_{0}^{t}%
\int_{0}^{t}|r -s |^{2H_{0}-2} \prod_{i=1}^d | B^i_{r}-B^i_{s } |
^{2H_{i}-2}\,dr \,ds \Biggr) \\
&<& \infty.\nonumber
\end{eqnarray}

\textit{Step} 2. We now prove that $u(t,x)$ is a weak
solution of (\ref{e.1.1}) in the sense of Definition \ref{def}.
Suppose that $\varphi$ is a smooth function with compact support. We
know that
%
%e4.7 ###
%
\begin{eqnarray} \label{a3}
&&
\int_{\mathbb{R}^d}u^{\varepsilon,\delta}(t,x)\varphi(x)\,dx \nonumber\\
&&\qquad=\int
_{%
\mathbb{R}^d}f(x)\varphi(x)\,dx+ \frac12 \int_{0}^{t}\int_{\mathbb{R}^d}u
^{\varepsilon,\delta} (t,x) \Delta\varphi(x)\,dx\,ds \\
&&\qquad\quad{}+\int_{0}^{t}\int_{\mathbb{R}^d}u ^{\varepsilon,\delta
}(t,x)\varphi(x)
\dot{W}^{\varepsilon,\delta}(s,x)\,ds\,dx.\nonumber
\end{eqnarray}
Therefore, it suffices to prove that
\[
\lim_{\varepsilon,\delta\downarrow0}\int_{0}^{t}\int_{\mathbb{R}^d}u
^{\varepsilon,\delta}(s,x)\varphi(x)\dot{W}^{\varepsilon,\delta
}(s,x)\,ds\,dx=\int_{0}^{t}\int_{\mathbb{R}^d}u(s,x)\varphi(x)W(ds,dx)
\]
in probability. From (\ref{a3}) and (\ref{a1}), it follows that $%
\int_{0}^{t}\int_{\mathbb{R}^d}u ^{\epsilon,\delta}(s,x)\varphi
(x)\dot{W}%
^{\epsilon,\delta}(s$,\break $x)\,ds\,dx$ converges in $L^{2}$ to the random
variable%
\[
G=\int_{\mathbb{R}^d}u(t,x)\varphi(x)\,dx-\int_{\mathbb
{R}^d}f(x)\varphi
(x)\,dx- \frac12 \int_{0}^{t}\int_{\mathbb{R}^d}u(t,x)\Delta\varphi(x)\,dx\,ds
\]
as $\varepsilon$ and $\delta$ tend to zero. Hence, if
\[
B_{\varepsilon,\delta}=\int_{0}^{t}\int_{\mathbb{R}^d}\bigl(u
^{\epsilon
,\delta}(s,x)-u(s,x)\bigr)\varphi(x)\dot{W}^{\varepsilon,\delta}(s,x)\,ds\,dx
\]
converges in $L^{2}$ to zero, then
\[
\int_{0}^{t}\int_{\mathbb{R}^d}u(s,x)\varphi(x)\dot W^{{\varepsilon
}, {
\delta}}\,ds\,dx = \int_{0}^{t}\int_{\mathbb{R}^d}u^{{\varepsilon},
{\delta}}
(s,x)\varphi(x)\dot W^{{\varepsilon}, {\delta}}\,ds\,dx -B_{{\varepsilon},
{%
\delta}}
\]
converges to $G$ in $L^2$. Thus, $u(s,x)\varphi(x)$ will be Stratonovich
integrable and we will have
\[
\int_{0}^{t}\int_{\mathbb{R}^d}u(s,x)\varphi(x)W(ds,dx)=G ,
\]
which will complete the proof. In order to show the convergence to zero
of $%
B_{ \varepsilon,\delta}$, we will express the product $(u
^{\varepsilon
,\delta}(s,x)-u(s,x))\dot{W}^{\varepsilon,\delta}(s,x)$ as the sum
of a
divergence integral plus a trace term [see (\ref{a5})]:
\begin{eqnarray*}
&& \bigl(u ^{\varepsilon,\delta}(s,x)-u(s,x)\bigr)\dot{W}^{\varepsilon,\delta}(s,x)
\\
&&\qquad = \int_{0}^{t}\int_{\mathbb{R}^d}\bigl(u ^{\varepsilon,\delta
}(s,x)-u(s,x)\bigr)\varphi_{\delta}(s-r)p_{\varepsilon}(x-z)\delta
W_{r,z} \\
&&\qquad\quad{} + \bigl\langle D\bigl(u ^{\varepsilon,\delta}(s,x)-u(s,x)\bigr),\varphi_{\delta
}(s-\cdot)p_{\varepsilon}(x-\cdot)\bigr\rangle_{\mathcal{H}} .
\end{eqnarray*}
We then have
%
%e4.8 ###
%
\begin{eqnarray} \label{BB}\hspace*{12pt}
B_{\varepsilon,\delta} &=& \int_{0}^{t}\int_{\mathbb{R}^d}\phi
_{r,z}^{\epsilon,\delta}\delta W_{r,z} \nonumber\\
&&{} +\int_{0}^{t}\int_{\mathbb{R}^d}\varphi(x)\bigl\langle D\bigl(u
^{\varepsilon
,\delta}(s,x)-u(s,x)\bigr),\varphi_{\delta}(s-\cdot)p_{\varepsilon
}(x-\cdot
)\bigr\rangle_{\mathcal{H}}\,ds\,dx \\
&=&B_{\varepsilon,\delta}^{1}+B_{\varepsilon,\delta}^{2},\nonumber
\end{eqnarray}
where
\[
\phi_{r,z}^{\varepsilon,\delta}=\int_{0}^{t}\int_{\mathbb{R}^d}\bigl(u
^{\varepsilon,\delta}(s,x)-u(s,x)\bigr)\varphi(x)\varphi_{\delta
}(s-r)p_{\varepsilon}(x-z)\,ds\,dx
\]
and $\delta(\phi^{\varepsilon,\delta})= \int_{0}^{t}\int_{\mathbb{R}
^d}\phi_{r,z}^{\varepsilon,\delta}\delta W_{r,z}$ denotes the divergence
or the Skorokhod integral of $\phi^{\varepsilon,\delta}$.

\textit{Step} 3. For the term $B_{\varepsilon
,\delta
}^{1}$, we use the following $L^{2}$ estimate for the Skorokhod integral:
%
%e4.9 ###
%
\begin{equation} \label{e.4.2}
E[(B_{\varepsilon,\delta}^{1})^{2}]\leq E(\Vert\phi^{\varepsilon
,\delta
}\Vert_{\mathcal{H}}^{2})+E(\Vert D\phi^{\varepsilon,\delta}\Vert_{
\mathcal{H}\otimes\mathcal{H}}^{2}) .
\end{equation}
The first term in (\ref{e.4.2}) is estimated as follows:%
%
%e4.10 ###
%
\begin{eqnarray}\label{a7}
&&E(\Vert\phi^{\varepsilon,\delta}\Vert_{\mathcal{H}}^{2})\nonumber\\
&&\qquad=
\int_{0}^{t}\int_{\mathbb{R}^d}\int_{0}^{t}\int_{\mathbb{R}^d}E \bigl[ \bigl(u
^{\varepsilon,\delta}(s,x)-u(s,x)\bigr)\nonumber\\
&&\qquad\quad\hspace*{75.4pt}{}\times\bigl(u ^{\varepsilon,\delta
}(r,y)-u(r,y)\bigr)%
\bigr]\varphi(x)\varphi(y) \\
&&\qquad\quad\hspace*{75.4pt}\hspace*{-14pt}{}\times\langle\varphi_{\delta}(s-\cdot)p_{\varepsilon}(x-\cdot
),\nonumber\\
&&\qquad\quad\hspace*{81.43pt}\varphi
_{\delta}(r-\cdot)p_{\varepsilon}(y-\cdot)\rangle_{\mathcal{H}}\,ds\,dx\,dr\,dy.\nonumber
\end{eqnarray}
Using Lemmas \ref{lemma2} and \ref{lemma3}, we can write%
%
%e4.11 ###
%
\begin{eqnarray} \label{a8}
&&\langle\varphi_{\delta}(s-\cdot)p_{\epsilon}(x-\cdot),\varphi
_{\delta}(r-\cdot)p_{\varepsilon}(y-\cdot)\rangle_{\mathcal{H}}
\nonumber
\\
&&\qquad=\alpha_H \biggl( \int_{[0,t]^{2}}\varphi_{\delta}(s-\sigma)\varphi
_{\delta}(r- \tau)|\sigma- \tau|^{2H_{0}-2}\,d \sigma \,d\tau\biggr)
\nonumber\\[-8pt]\\[-8pt]
&&\qquad\quad{}\times\Biggl( \int_{\mathbb{R}^{2d}}p_{\varepsilon}(x-z )p_{\varepsilon
}(y-w) \prod_{i=1}^d |z_i -w_i|^{2H_{i}-2}\,dz \,dw \Biggr) \nonumber\\
&&\qquad\leq C|s-r|^{2H_{0}-2} \prod_{i=1}^d |x-y|^{2H_{i}-2}\nonumber
\end{eqnarray}
for some constant $C>0$. As a consequence, the integrand on the right-hand
side of (\ref{a7}) converges to zero as $\varepsilon$ and
$\delta$
tend to zero for any $s$, $r$, $x$, $y$ due to (\ref{a1}). From (\ref{a2a}),
we get%
%
%e4.12 ###
%
\begin{eqnarray}\label{a9}
&&\sup_{\varepsilon,\delta}\sup_{x\in\mathbb{R}^d}\sup_{0\leq
s\leq
t}E ( u ^{\epsilon,\delta}(s,x) ) ^{2}\nonumber\\[-8pt]\\[-8pt]
&&\qquad\leq\| f \|
_{\infty}^{2}\sup_{\varepsilon,\delta}\sup_{x\in\mathbb
{R}^d}\sup
_{0\leq
s\leq t}E\exp( 2V_{s,x}^{\varepsilon,\delta} ) <\infty.\nonumber
\end{eqnarray}
Hence, from (\ref{a8}) and (\ref{a9}), we get that the integrand on the
right-hand side of (\ref{a7}) is bounded by $C|s-r|^{2H_{0}-2}
\prod_{i=1}^d |x_i-y_i|^{2H_{i}-2}$ for some constant $C>0$.
Therefore, by
dominated convergence, we get that $E(\Vert\phi^{\varepsilon,\delta
}\Vert
_{\mathcal{H}}^{2})$ converges to zero as $\varepsilon$ and $\delta$ tend
to zero.

\textit{Step} 4. On the other hand, we have
\[
D(u ^{\varepsilon,\delta}(t,x))=E^{B} [ f(B_{t}+x)\exp
(V_{t,x}^{\varepsilon,\delta})A_{t,x}^{\varepsilon,\delta} ] ,
\]
where $A_{t,x}^{\varepsilon,\delta}$ is defined as in (\ref{e2}).
Therefore,%
%
%e4.13 ###
%
\begin{eqnarray} \label{e.4.12a}
&&E\langle D(u ^{\varepsilon,\delta}(t,x)),D(u ^{\varepsilon^{\prime
},\delta^{\prime}}(t,x))\rangle_{\mathcal{H}} \nonumber\\
&&\qquad=E^{W}E^{B} \bigl( f(B_{t}^{1}+x)f(B_{t}^{2}+x) \exp\bigl(V_{t,x}^{\varepsilon,\delta
}(B^{1})+V_{t,x}^{\varepsilon
,\delta}(B^{2})\bigr)\\
&&\qquad\quad\hspace*{133.9pt}{} \times\langle A_{t,x}^{\varepsilon,\delta
}(B^{1}),A_{t,x}^{\varepsilon^{\prime},\delta^{\prime
}}(B^{2})\rangle_{\mathcal{H}} \bigr),\nonumber
\end{eqnarray}
where $B^{1}$ and $B^{2}$ are two independent $d$-dimensional Brownian
motions and where $E^B$ denotes the expectation with respect to $(B^1,
B^2)$%
. Then, from the previous results it is easy to show that%
%
%e4.14 ###
%
\begin{eqnarray} \label{b1}
&&\lim_{\varepsilon,\delta\downarrow0}E\langle D(u ^{\varepsilon
,\delta
}(t,x)),D(u ^{\varepsilon^{\prime},\delta^{\prime}}(t,x))\rangle_{%
\mathcal{H}} \nonumber\\
&&\qquad=E \Biggl[ f(B_{t}^{1}+x)f(B_{t}^{2}+x)\nonumber\\[-8pt]\\[-8pt]
&&\qquad\quad\hspace*{12pt}{}\times\exp\Biggl( \frac{\alpha_H}{2} \sum_{j,k=1}^{2}\int_{0}^{t}%
\int_{0}^{t}|s-r|^{2H_{0}-2} \prod_{i=1}^d
|B_{s}^{j,i}-B_{r}^{k,i}|^{2H_{i}-2}\,ds\,dr \Biggr) \nonumber\\
&&\qquad\quad\hspace*{66.8pt}{}\times\alpha_H \int_{0}^{t}\int_{0}^{t}|s-r|^{2H_{0}-2}
\prod_{i=1}^d |B_{s}^{1,i}-B_{r}^{2,i}|^{2H_{i}-2}\,ds\,dr \Biggr]
.\hspace*{-12pt}\nonumber
\end{eqnarray}
This implies that $u ^{\varepsilon,\delta}(t,x)$ converges in the
space $%
\mathbb{D}^{1,2}$ to $u(t,x)$ as $\delta\downarrow0$ and
$\varepsilon
\downarrow0$. Letting ${\varepsilon}^{\prime}={\varepsilon}$ and
${\delta}%
^{\prime}={\delta}$ in (\ref{e.4.12a}) and using the same
argument as
for (\ref{a9}), we obtain
\[
\sup_{\varepsilon,\delta}\sup_{x\in\mathbb{R}^d}\sup_{0\leq
s\leq
t}E \| D(u ^{\varepsilon,\delta}(s,x)) \| _{\mathcal{H}%
}^{2}<\infty.
\]
Then,
\begin{eqnarray*}
&&E\Vert D\phi^{\varepsilon,\delta}\Vert_{\mathcal{H}\otimes
\mathcal
{H}%
}^{2} \\
&&\qquad=\int_{0}^{t}\int_{\mathbb{R}^d}\int_{0}^{t}\int_{\mathbb
{R}%
}E\bigl\langle D\bigl(u ^{\varepsilon,\delta}(s,x)-u(s,x)\bigr),D\bigl(u ^{\varepsilon
,\delta
}(r,y)-u(r,y)\bigr)\bigr\rangle_{\mathcal{H}} \\
&&\qquad\quad\hspace*{58.6pt}{}\times\varphi(x)\varphi(y)\langle\varphi_{\delta}(s-\cdot
)p_{\varepsilon}(x-\cdot),\\
&&\qquad\quad\hspace*{120.9pt}\varphi_{\delta}(r-\cdot)p_{\varepsilon
}(y-\cdot)\rangle_{\mathcal{H}}\,ds\,dx\,dr\,dy
\end{eqnarray*}
converges to zero as $\varepsilon$ and $\delta$ tend to zero. Hence,
by (%
\ref{e.4.2}), $B_{\varepsilon,\delta}^{1}$ converges to zero in $L^{2}$
as $\varepsilon$ and $\delta$ tend to zero.

\textit{Step} 5. The second summand in the right-hand
side of (\ref{BB}) can be written as%
\begin{eqnarray*}
B_{\varepsilon,\delta}^{2} &=&\int_{0}^{t}\int_{\mathbb
{R}^d}\varphi
(x)\bigl\langle D\bigl(u ^{\varepsilon,\delta}(s,x)-u(s,x)\bigr),\varphi_{\delta
}(s-\cdot)p_{\varepsilon}(x-\cdot)\bigr\rangle_{\mathcal{H}}\,ds\,dx \\
&=&\int_{0}^{t}\int_{\mathbb{R}^d}\varphi(x)E^{B} \bigl(
f(B_{s}^{x})\exp
( V_{s,x}^{\varepsilon,\delta} ) \langle A_{s,x}^{\varepsilon
,\delta},\varphi_{\delta}(s-\cdot)p_{\epsilon}(x-\cdot)\rangle_{%
\mathcal{H}} \bigr) \,ds\,dx \\
&&{}-\int_{0}^{t}\int_{\mathbb{R}}\varphi(x)E^{B} \bigl( f(B_{s}^{x})\exp
( V_{s,x} ) \\
&&\hspace*{79.3pt}{}\times\langle\delta(B_{s-\cdot}^{x}-\cdot),\varphi
_{\delta}(s-\cdot)p_{\varepsilon}(x-\cdot)\rangle_{\mathcal{H}} \bigr)
\,ds\,dx \\
&=&B_{\varepsilon,\delta}^{3}-B_{\varepsilon,\delta}^{4},
\end{eqnarray*}
where%
\begin{eqnarray*}
&&\langle A_{s,x}^{\varepsilon,\delta},\varphi_{\delta}(s-\cdot
)p_{\epsilon}(x-\cdot)\rangle_{\mathcal{H}}\\
&&\qquad= \alpha
_H\int_{[0,s]^{3}}\int_{\mathbb{R}^{2d}} |r -v|^{2H_{0}-2}\\
&&\qquad\quad\hspace*{62.8pt}{}\times \prod_{i=1}^d
|y_i-z_i|^{2H_{i}-2} \varphi_{\delta}(s-r )p_{\varepsilon}(B_{r}^x-y ) \\
&&\qquad\quad\hspace*{62.8pt}{}\times\varphi_{\delta}(s-v)p_{\varepsilon}(x-z)\,dy \,dz\,dr\,dr\,dv
\end{eqnarray*}
and%
\begin{eqnarray*}
&&\langle\delta(B_{s-\cdot}^{x}-\cdot),\varphi_{\delta}(s-\cdot
)p_{\varepsilon}(x-\cdot)\rangle_{\mathcal{H}} \\
&&\qquad=\alpha_H \int_{[0,s]^{2}}\int_{\mathbb{R}^d} v^{2H_{0}-2} \prod_{i=1}^d
|B_{r} ^{x_i} -y_i|^{2H_{i}-2}\varphi_{\delta}(r-v)p_{\varepsilon
}(x-y)\,dy\,dv \,dr.
\end{eqnarray*}
Lemma \ref{lemma2} and Lemma \ref{lemma3} imply that%
%
%e4.15 ###
%
\begin{equation} \label{c1}
\langle A_{s,x}^{\varepsilon,\delta},\varphi_{\delta}(s-\cdot
)p_{\varepsilon}(x-\cdot)\rangle_{\mathcal{H}}\leq
C\int_{0}^{s}r^{2H_{0}-2} \prod_{i=1}^d |B^i_{r}|^{2H_{i}-2}\,dr
\end{equation}
and%
%
%e4.16 ###
%
\begin{equation} \label{c2}\quad
\langle\delta(B_{s-\cdot}^{x}-\cdot),\varphi_{\delta}(s-\cdot
)p_{\varepsilon}(x-\cdot)\rangle_{\mathcal{H}}\leq
C\int_{0}^{s}r^{2H_{0}-2} \prod_{i=1}^d |B^i_{r}|^{2H_{i}-2}\,dr
\end{equation}
for some constant $C>0$. Then, from (\ref{c1}), (\ref{c2}) and the
fact that the random variable $\int_{0}^{s}r^{2H_{0}-2} \prod_{i=1}^d
|B^i_{r}|^{2H_{i}-2}\,dr$ is square integrable because of Lemma \ref{lemma4},
we can apply the dominated convergence theorem and get that
$B_{\varepsilon
,\delta}^{3}$ and $B_{\varepsilon,\delta}^{4}$ both converge in
$L^{2}$ to%
\[
\alpha_H\int_{0}^{t}\int_{\mathbb{R}^d}\varphi(x)E^{B} \Biggl(
f(B_{s}^{x})\exp( V_{s,x} ) \int_{0}^{s}r^{2H_{0}-2}
\prod_{i=1}^d |B^i_{r}|^{2H_{i}-2}\,dr \Biggr) \,ds\,dx
\]
as $\varepsilon$ and $\delta$ tend to zero. Therefore,
$B_{\varepsilon
,\delta}^{2}$ converges to zero in $L^{2}$ as $\varepsilon$ and
$\delta$
tend to zero. This completes the proof.
\end{pf}

We can also show that the process $u(t,x)$ given in (\ref{utx}) is a
mild solution to~(\ref{e.1.1}), in the sense that the
following equation holds:
\[
u(t,x)= p_tf(x)+ \int_0^t \int_{\mathbb{R}^d} p_{t-s}(x-y) u(s,y) \,dW_{s,y},
\]
where $p_t$ denotes the heat kernel and $ p_tf(x)=\int_{\mathbb{R}^d}
p_{t }(x-y) f(y)\,dy$. In fact, as in the proof of Theorem \ref{th1}, we
need to show that
\[
\int_0^t \int_{\mathbb{R}^d} p_{t-s}(x-y) \bigl(u(s,y)-u^{\varepsilon
,\delta
}(s,y)\bigr) \,dW^{\varepsilon,\delta}_{s,y}
\]
converges to zero in $L^2$. This can be proven with the same arguments
as in the proof of Theorem \ref{th1}, replacing $\varphi$ by the heat
kernel. For instance, instead of the estimate (\ref{a8}), we should have
\begin{eqnarray*}
&& \int_0^t \int_0^t \int_{\mathbb{R}^{2d}} p_{t-r}(x-y)p_{t-s}(x-z )
|s-r|^{2H_{0}-2} \prod_{i=1}^d |y-z|^{2H_{i}-2} \,dy\,dz\,dr\,ds \\
&&\qquad = \int_0^t \int_0^t |s-r|^{2H_{0}-2} E \Biggl( \prod_{i=1}^d |
B^{1,i}_{t-r} -B^{2,i}_{t-s} |^{2H_{i}-2} \Biggr)\,dr\,ds<\infty.
\end{eqnarray*}
We omit the details of this proof.
\begin{remark}
The uniqueness of the solution remains to be investigated in a
future work. The definition of the Stratonovich integral as a limit in
probability makes the uniqueness problem nontrivial and it is not clear
how to proceed.
\end{remark}

As a corollary of Theorem \ref{th1}, we obtain the following result.
\begin{corollary}
\label{t.4.4} Suppose that $2H_{0}+\sum_{i=1}^{d}H_{i}>d+1$. Then, the
solution $%
u(t,x)$ given by (\ref{utx}) has finite moments of all orders.
Moreover, for any positive integer $p$, we have
%
%e4.17 ###
%
\begin{eqnarray}\label{e.4.7}\qquad
&&E ( u(t,x)^{p} ) \nonumber\\
&&\qquad= E \Biggl(\prod_{j=1}^{p}f(B_{t}^{j}+x)\nonumber\\[-8pt]\\[-8pt]
&&\qquad\quad\hspace*{29.4pt}{}\times\exp\Biggl[ \frac{\alpha_{H}}{2}\sum_{j,k=1}^{p}\int_{0}^{t}%
\int_{0}^{t}|s-r|^{2H_{0}-2}\nonumber\\
&&\qquad\quad\hspace*{133pt}{}\times
\prod_{i=1}^{d}|B_{s}^{j,i}-B_{r}^{k,i}|^{2H_{i}-2}\,ds\,dr \Biggr] \Biggr) ,
\nonumber
\end{eqnarray}
where $B_{1},\ldots,B_{p}$ are independent $d$-dimensional standard
Brownian motions.
\end{corollary}
\begin{remark}
In the previous work \cite{hunualart}, a formula similar
to (\ref{e.4.7}) was obtained in
the special case $H_1=\cdots=H_d=\frac12$, without the condition $%
2H_{0}+\sum_{i=1}^{d}H_{i}>d+1$. This type of formula was proven
assuming $d=1$ and $H_0>\frac34$. In the case of the Skorokhod-type
equation, a formula for the moments of the solution similar to
(\ref{e.4.7}) was established in \cite{hunualart} if
$d=1$ or $2$, $H_0>\frac12$ and $t$ is small enough.
\end{remark}

%s5 ###
\section{Behavior of the Feynman--Kac formula}\label{sec5}

In this section, we present two applications of the
Feynman--Kac formula.

%s5.1 ###
\subsection{\texorpdfstring{H\"{o}lder continuity of the solution}{Holder continuity of the solution}}

In this subsection, we study the H\"older continuity of the solution of
(\ref{e.1.1}). The main result of this section is the following
theorem.
\begin{theorem}
\label{t.6.1} Suppose that $2H_0+\sum_{i=1}^d H_{i}>d+1$ and let
$u(t,x)$ be
the solution of (\ref{e.1.1}). Then, $u(t,x)$ has a continuous
modification such that for any $\rho\in(0, \frac\kappa2 )$
[where $\kappa$ is defined as in (\ref{kappa})] and any compact
rectangle $I \subset\mathbb{R}_+\times\mathbb{R}^d$, there exists a
positive random variable $K_I $ such that almost surely, for any $%
(s,x),(t,y)\in I$, we have
\[
| u(t,y )- u(s,x )|\le K_I (|t-s|^\rho+|y-x|^{2\rho}).
\]
\end{theorem}
\begin{pf}
The proof involves several steps.

\textit{Step} 1. Recall that $V_{t,x}=\int_0^t\int
_{%
\mathbb{R}^d} \delta(B_{t-r}^x-y)W(dr,dy)$ denotes the random variable
introduced in (\ref{e.3.4}) and
\[
u(t,x)= E^B (f(B^x_t) \exp(V_{t,x}).
\]
Set $V =V_{s,x}$ and $\tilde V=V _{t,y }$. We can then write
\begin{eqnarray*}
&& E^W|u(s,x )-u(t,y)|^p \\
&&\qquad = E^W|E^B (e^{V }-e^{\tilde V})|^p \\
&&\qquad \le E^W \bigl(E^B \bigl[|\tilde V -V |e^{ \max(V , \tilde V )}\bigr] \bigr)^p \\
&&\qquad \le E^W \bigl[ \bigl(E^Be^{2 \max(V , \tilde V ) } \bigr)^{p/2} \bigl(%
E^B(\tilde V -V )^2 \bigr)^{p/2} \bigr] \\
&&\qquad \le\bigl[E^W E^B e^{2p \max(V , \tilde V ) } \bigr]^{1/2} \bigl[E^W%
\bigl(E^B (\tilde V -V )^2 \bigr)^p \bigr]^{1/2}.
\end{eqnarray*}
Applying Minkowski's inequality, the equivalence between the $L^2$-norm and
the $L^p$-norm for a Gaussian random variable and using the exponential
integrability property (\ref{e21}), we obtain
%
%e5.1 ###
%
\begin{eqnarray} \label{j1}
E^W|u(s,x)-u(t,y)|^p &\le& C \bigl[E^W \bigl(E^B (\tilde V -V )^2 \bigr)^p \bigr]^%
{1/2} \nonumber\\[-8pt]\\[-8pt]
&\le& C_p [E^B E^W|\tilde V -V |^2 ]^{p/2}.\nonumber
\end{eqnarray}
In a similar way to (\ref{e.2.14}), we can deduce the following
formula for
the conditional variance of $\tilde V -V$:
%
%e5.2 ###
%
\begin{eqnarray} \label{j2}
&&E^W|\tilde V -V |^2 \nonumber\\
&&\qquad= \alpha_H E^B \Biggl( \int_0^{s}\int_0^{s} |r
-v|^{2H_0-2} \prod_{i=1}^d |B^{i}_{s-r }-B^{i}_{s-v}|^{2H_i-2}\,dr\,dv
\nonumber\\
&&\qquad\quad\hspace*{34.8pt}{} + \int_0^{t }\int_0^{t } |r-v|^{2H_0-2} \prod_{i=1}^d|B^{i}_{t -r
}-B^{i}_{t -v}|^{2H_i-2}\,dr\,dv \nonumber\\[-8pt]\\[-8pt]
&&\qquad\quad\hspace*{34.8pt}{} -2\int_0^{s}\int_0^{t } |r-v|^{2H_0-2}\nonumber\\
&&\qquad\quad\hspace*{83.5pt}{}\times \prod_{i=1}^d|B_{s-r
}^{i}-B_{t-v}^{i} +x_i -y_i|^{2H_i-2}\,dr\,dv \Biggr) \nonumber\\
&&\qquad:= \alpha_H C(s, t , x , y).\nonumber
\end{eqnarray}

\textit{Step} 2. $ $ Fix $1\leq j\leq d$. Let us
estimate $C(s,t,x,y)$ when $s=t$ and $x_{i}=y_{i}$ for all $i\not=j$. We
can write
%
%e5.3 ###
%
\begin{eqnarray} \label{eqa1}
&&C(t,t,x,y)\nonumber\\
&&\qquad=2\int_{0}^{t}\int_{0}^{t}|r-v|^{\kappa-1}\\
&&\qquad\quad\hspace*{34.8pt}{}\times\prod_{i\not
=j}^{d}E%
(|\xi|^{2H_{i}-2} )E (|\xi|^{2H_{j}-2}-|z+\xi|^{2H_{j}-2} )%
\,dr\,dv,\nonumber
\end{eqnarray}
where $z=\frac{x_{j}-y_{j}}{\sqrt{|r-v|}}$ and $\xi$ is a standard normal
variable. Set $\beta_{j}=2H_{j}+1>2$. By Lemma \ref{Lem7.6}, the
factor $E%
(|\xi|^{2H_{j}-2}-|z+\xi|^{2H_{j}-2} )$ can be bounded by a
constant if $|r-v|\leq(x_{j}-y_{j})^{2}$ and it can be bounded by $%
C|x_{j}-y_{j}|^{\beta_{j}}|r-v|^{-\beta_{j}/2}$ if $%
|r-s|>(x_{j}-y_{j})^{2} $. In this way, we obtain
\begin{eqnarray*}
C(t,t,x,y)&\leq& C\int_{\{0<r,v<t,|r-v|\leq
(x_{j}-y_{j})^{2}\}}|r-v|^{\kappa-1}\,dr\,dv \\
&&{} +C|x_{j}-y_{j}|^{\beta
_{j}}\int_{\{0<r,v<t,|r-v|>(x_{j}-y_{j})^{2}\}}|r-v|^{\kappa-1-\beta
_{j}/2}\,dr\,dv \\
&\leq& C|x_{j}-y_{j}|^{2\kappa}.
\end{eqnarray*}
So, from (\ref{j1}), we have
%
%e5.4 ###
%
\begin{equation} \label{e.6.1}
E^{W}|u(t,x)-u(t,y)|^{p}\leq C|x_{j}-y_{j}|^{\kappa p}.
\end{equation}

\textit{Step} 3. Now, suppose that $s<t $ and $x =y$.
Set $\delta=\sum_{i=1}^d H_i-d$. We have
\begin{eqnarray*}
&&C(s,t,x,x) \\
&&\qquad=C \biggl[ \int_{s} ^{t } \int_{s} ^{t } |r-v |^{\kappa-1 }\,dr\,dv \\
&&\qquad\quad\hspace*{12pt}{}+\int_0^{s}\int_0^{t }|r-v|^{2H_0-2} (|r -v|^{\delta}-|r -v+
t-s|^{\delta} )\,dr\,dv \biggr] .
\end{eqnarray*}
The first integral is $O((t-s)^{\kappa+1})$ when $t-s$ is small. For the
second integral, we use the change of variable $\sigma=r -v$, $v=\tau
$ and
we have
\begin{eqnarray*}
&&\int_0^{s}\int_0^{t }|r-v|^{2H_0-2} (|r -v|^{\delta}-|r -v+
t-s|^{\delta
} )\,dr\,dv \\
&&\qquad\le \int_0^{t } d\tau\int_{-t }^{s} |\sigma|^{2H_0-2} \bigl| |\sigma
|^{\delta}-|\sigma+ t-s |^{\delta} \bigr|\,d\sigma\\
&&\qquad= t \biggl[ \int_0^{s} \sigma^{2H_0-2} \bigl(\sigma^{\delta}-(\sigma+
t-s)^{\delta} \bigr)\,d\sigma\\
&&\qquad\quad\hspace*{8pt}{}+\int_{-t}^{s-t} (-\sigma)^{2H_0-2} \bigl((-\sigma-
t+s)^{\delta}-(-\sigma)^{\delta} \bigr)\,d\sigma\\
&&\qquad\quad\hspace*{23.3pt}{}+\int_{s-t}^0 (-\sigma)^{2H_0-2} |(-\sigma
)^{\delta}-(\sigma+t-s)^{\delta} |\,d\sigma\biggr] \\
&&\qquad= t [A^{\prime}+B^{\prime}+C^{\prime}].
\end{eqnarray*}
For the first term in the above decomposition, we can write
\begin{eqnarray*}
A^{\prime}&=&(t-s)^{\kappa-1} \int_0^{{t_1}/({t-s})} \sigma^{2H_0-2}
\bigl(%
\sigma^{\delta}-(\sigma+1)^{\delta} \bigr)\,d\sigma\\
&\le& (t-s)^{\kappa-1} \int_0^\infty\sigma^{2H_0-2} \bigl(%
\sigma^{\delta}-(\sigma+1)^{\delta} \bigr)\,d\sigma\\
&\le& C(t-s)^{\kappa},
\end{eqnarray*}
because $2H_0 + \sum_{i=1}^d -d-3<-1$. Similarly, we can get that
\[
B^{\prime}\le(t-s)^{\kappa} \int_1^\infty\sigma^{2H_0-2} \bigl(\sigma
^{\delta}-(\sigma+1)^{\delta} \bigr)\,d\sigma.
\]
Finally,
\[
C^{\prime}\le\int_0^{t-s} \sigma^{2H_0-2} \bigl(\sigma^{\delta}+(
t-s-\sigma)^{\delta} \bigr)\,d\sigma=C (t-s)^{\kappa}.
\]
So, we have
%
%e5.5 ###
%
\begin{equation} \label{e.6.2}
E^W|u(s,x)-u(t,y)|^p\le C (t-s)^{\kappa/2 p}.
\end{equation}

\textit{Step} 4. Combining equations (\ref{e.6.1}) and
(\ref{e.6.2}) with the estimates (\ref{j1}) and (\ref{j2}), the
result of this theorem can now be concluded from Theorem 1.4.1 of
Kunita %
\cite{kunita} if we choose $p$ large enough.
\end{pf}

%It is well-known that when $d=1$ and $H_0=H_1=\frac12$, the solution
%of the equation in
%the It\^o sense is jointly
%H\"older continuous with exponent ${\alpha}$ in time variable and $2{
%is space variable for any ${\alpha}\in(0, 1/4)$. It is obvious that
%when
%Theorem \ref{t.6.1} is formally applied to this case, one cannot
%obtain the
%known result. This makes us believe that the modulus of continuity may
%be
%improved. The technique in \cite{VV} may be combined with the explicit
%expression (\ref{utx}).

%s5.2 ###
\subsection{Regularity of the density}

In this subsection, we shall use the Feynman--Kac formula established
in the
previous section to show that for any $t$ and $x$, the probability law of
the solution $u(t,x)$ of (\ref{e.1.1}) has a smooth density with
respect to the Lebesgue measure. To this end, we shall show that $ \|
Du(t,x) \| _{\mathcal{H}}$ has negative moments of all orders.
\begin{theorem} \label{th2}
Suppose that $2H_0+\sum_{i=1}^d H_{i}>d+1$. Fix $t>0$ and $x\in
\mathbb
{R}^d$%
. Assume that for any positive number $p$, $E|f(B_{t}+x)|^{-p}<\infty$.
Then, the law of $u(t,x)$ has a smooth density.
\end{theorem}
\begin{pf}
From Theorem \ref{th1}, we can write
\[
u(t,x)=E^{B} [ f(B_{t}^{x})\exp( V_{t,x} ) ] .
\]
The Malliavin derivative of the solution is given by
\[
D_{r,y}u(t,x)=E^{B} [ f(B_{t}^{x})\exp( V_{t,x} ) \delta
(B_{t-r}^{x}-y) ] .
\]
It is not difficult to show that $u(t,x)\in\mathbb{D}^{\infty}$.
Thus, by
the general criterion for the smoothness of densities (see \cite{nualart}),
it suffices to show that $E ( \| Du(t,x) \| _{\mathcal{H}%
}^{-2p} ) <\infty$ for any $t>0$ and $x\in\mathbb{R}^{d}$. We have
\begin{eqnarray*}
\| Du(t,x) \| _{\mathcal{H}}^{2} &=&E^{B} \bigl[
f(B_{t}^{1}+x)f(B_{t}^{2}+x)\exp\bigl( V_{t,x}(B^{1})+V_{t,x}(B^{2}) \bigr)
\\
&&{}\hspace*{89pt}\times\langle\delta(B_{t-r}^{1,x}-y) ,\delta
(B_{t-r}^{2,x}-y)\rangle_{\mathcal{H}} \bigr] \\
&=&\alpha_{H}E^{B} \Biggl[ f(B_{t}^{1}+x)f(B_{t}^{2}+x)\exp\bigl(
V_{t,x}(B^{1})+V_{t,x}(B^{2}) \bigr) \\
&&\hspace*{33.1pt}{}\times
\int_{0}^{t}\int_{0}^{t}|r-s|^{2H_{0}-2}%
\prod_{i=1}^{d}|B_{t-r}^{1,i}-B_{t-s}^{2,i}|^{2H_{i}-2}\,dr\,ds \Biggr] ,
\end{eqnarray*}
where $B^{1}$ and $B^{2}$ are independent $d$-dimensional Brownian motions.
By Jensen's inequality, we have, for any $p>0$,
that
\begin{eqnarray*}
&& \| Du(t,x) \| _{\mathcal{H}}^{-2p} \\
&&\qquad\leq ( \alpha_{H} ) ^{-p}E^{B} \Biggl[
|f(B_{t}^{1}+x)f(B_{t}^{2}+x)|^{-p}\exp\bigl( -p [
V_{t,x}(B^{1})+V_{t,x}(B^{2}) ] \bigr) \\
&&\qquad\quad\hspace*{56.3pt}{}\times\biggl(
\int_{0}^{t}\int_{0}^{t}|r-s|^{2H_{0}-2}%
\prod_{i=1}^{d}|B_{t-r}^{i,1}-B_{t-s}^{2,i}|^{2H_{i}-2}\,dr\,ds \biggr) ^{-p}\Biggr] .
\end{eqnarray*}
Hence, by H\"{o}lder's inequality, we obtain
\begin{eqnarray*}
&&E \| Du(t,x) \| _{\mathcal{H}}^{-2p} \\
&&\qquad\leq ( \alpha_{H} ) ^{-p} \bigl(
E|f(B_{t}^{1}+x)f(B_{t}^{2}+x)|^{-pp_{1}} \bigr) ^{{1}/{p_{1}}} \\
&&\qquad\quad{}\times\bigl( E\exp\bigl( -pp_{2} [ V_{t,x}(B^{1})+V_{t,x}(B^{2})%
] \bigr) \bigr) ^{{1}/{p_{2}}} \\
&&\hspace*{14.2pt}\qquad\quad{} \times\Biggl( E \Biggl(
\int_{0}^{t}\int_{0}^{t}|r-s|^{2H_{0}-2}%
\prod_{i=1}^{d}|B_{t-r}^{1,i}-B_{t-s}^{2,i}|^{2H_{i}-2}\,dr\,ds \Biggr)
^{-pp_{3}} \Biggr) ^{{1}/{p_{3}}} \\
&&\qquad=I_{1}I_{2}I_{3} ,
\end{eqnarray*}
where $\frac{1}{p_{1}}+\frac{1}{p_{2}}+\frac{1}{p_{3}}=1$. The first
factor, $%
I_{1}$, is finite by the assumption on $f$ and H\"{o}lder's inequality. The
second factor is finite by Theorem \ref{teo2}. Finally, from Jensen's
inequality, we have
\begin{eqnarray*}
I_{3}^{p_{3}} &=&E \Biggl[ t^{{-2pp_{3}}} \Biggl\{ \frac{1}{t^{2}}%
\int_{0}^{t}\int_{0}^{t}|r-s|^{2H_{0}-2}%
\prod_{i=1}^{d}|B_{t-r}^{1,i}-B_{t-s}^{2,i}|^{2H_{i}-2}\,dr\,ds \Biggr\}
^{-pp_{3}} \Biggr] \\
&\leq&E \Biggl[ t^{{-2pp_{3}-2}} \Biggl\{
\int_{0}^{t}\int_{0}^{t}|r-s|^{-(2H_{0}-2)pp_{3}}\\
&&\hspace*{88.9pt}{}\times\prod_{i=1}^{d}|B_{t-r}^{1,i}-B_{t-s}^{2,i}|^{-(2H_{i}-2)pp_{3}}\,dr\,ds \Biggr\}
\Biggr] \\
&\leq&C\int_{0}^{t}\int_{0}^{t}|r-s|^{-(2H_{0}-2)pp_{3}}E \Biggl\{
\prod_{i=1}^{d}|B_{t-r}^{1,i}-B_{t-s}^{2,i}|^{-(2H_{i}-2)pp_{3}} \Biggr\} \,dr\,ds
\\
&<&\infty.
\end{eqnarray*}
This completes the proof.
\end{pf}

%s6 ###
\section{The case $H_0>\frac{3}{4}, H_1=\frac{1}{2}$ and $d=1$}\label{sec6}

%s6.1 ###
\subsection{Preliminaries}

In this case, all the setup is the same as before, except that if $\phi
$ and
$\psi$ are functions in $\mathcal{E}$, then%
\begin{eqnarray*}
E ( W(\phi)W(\psi) ) &=& \langle\phi,\psi\rangle_{\mathcal{H}%
}\\
&=&\alpha_{H_{0}}
\int_{0}^{\infty}\int_{0}^{\infty}\int_{\mathbb{R}}\phi
(s,x)\psi(t,x)|s-t|^{2H_{0}-2}\,ds\,dt\,dx,
\end{eqnarray*}
where $\alpha_{H_{0}}=H_{0}(2H_{0}-1)$.

%s6.2 ###
\subsection{Definition and exponential integrability
of the stochastic Feynman--Kac functional}

Similarly, we also have the following theorem.
\begin{theorem}\label{theo61}
Suppose that $H_1=1/2 $ and $H_{0}>3/4$. Then, for any $\varepsilon>0 $
and $\delta>0$, $A_{t,x}^{\varepsilon,\delta}$ defined in (\ref{e2})
belongs to $\mathcal{H}$ and the family of random variables $%
V_{t,x}^{\epsilon,\delta}$ defined in (\ref{e3}) converges in $L^{2}$
to a
limit denoted by
%
%e6.1 ###
%
\begin{equation}\label{e.8.8}
V_{t,x}=\int_{0}^{t}\int_{\mathbb{R}}\delta(B_{t-r}^{x}-y)W(dr,dy) .
\end{equation}
Conditional on $B$, $V_{t,x}$ is a Gaussian random variable with mean $0$
and variance%
%
%e6.2 ###
%
\begin{equation} \label{e.8.9}
\operatorname{Var}^W(V_{t,x}) =\alpha_{H_{0}}\int_{0}^{t}\int_{0}^{t}|r
-s|^{2H_{0}-2}\delta( B_{r }-B_{s} )\,dr \,ds .
\end{equation}
\end{theorem}
\begin{pf}
Fix $\varepsilon$, $\varepsilon^{\prime}$, $\delta$ and $\delta
^{\prime
}>0$.
\begin{eqnarray*}
&&E^BE^W (V_{t,x}^{\epsilon,\delta},V_{t,x}^{\epsilon^{\prime},\delta
^{%
\prime}} ) \\
&&\qquad= E^B \langle A_{t,x}^{\varepsilon,\delta
},A_{t,x}^{\varepsilon^{\prime},\delta^{\prime}} \rangle_{\mathcal{
H}} \\
&&\qquad= \alpha_{H_{0}}E^B \biggl(\int_{[0,t]^{4}}\int_{\mathbb{R}%
}p_{\epsilon}(B_{s}^{x}-y)p_{\epsilon^{\prime}}(B_{r}^{x}-y)\varphi_{\delta}(t-s -u ) \\
&&\qquad\quad\hspace*{76.2pt}{} \times\varphi_{\delta^{\prime
}}(t-r-v) |u-v|^{2H_{0}-2}\,dy\,du\,dv\,ds\,dr \biggr) \\
&&\qquad= \alpha_{H_{0}} \biggl(\int_{[0,t]^{4}}E^Bp_{\epsilon
+\epsilon^{\prime}}(B_{s}-B_{r})\varphi_{\delta}(t-s-u) \\
&&\qquad\quad\hspace*{47.3pt}{} \times\varphi_{\delta^{\prime
}}(t-r-v) |u-v|^{2H_{0}-2}\,du\,dv\,ds\,dr \biggr) \\
&&\qquad= \alpha_{H_{0}} \biggl(\int_{[0,t]^{4}}\frac{1}{\sqrt{2\pi}}%
(\epsilon+\epsilon^{\prime}+|s-r| )^{-1/2}\varphi_{\delta}(t-s-u) \\
&&\qquad\quad\hspace*{48.1pt}{} \times\varphi_{\delta^{\prime}}(t-r-v)
|u-v|^{2H_{0}-2}\,du\,dv\,ds\,dr \biggr).
\end{eqnarray*}
By Lemma \ref{lemma3},
\begin{eqnarray*}
&&\int_{[0,t]^2} (\epsilon+\epsilon^{\prime}+|s-r| )^{-1/2} \times
\varphi_{\delta}(t-s-u)\varphi_{\delta^{\prime}}(t-r-v)
|u-v|^{2H_{0}-2}\,du\,dv \\
&&\qquad \le C |s-r |^{2H_0- 5/2}.
\end{eqnarray*}
Then, by the dominated convergence theorem, $E^BE^W (V_{t,x}^{\epsilon,
\delta},V_{t,x}^{\epsilon^{\prime},\delta^{\prime}} )$ converges to
\[
\frac{\alpha_{H_{0}}}{\sqrt{2\pi}} \int_{[0,t]^{2}} |s-r|^{2H_{0}-5/2}\,ds\,dr
\]
as $\varepsilon$, $\varepsilon^{\prime}$, $\delta$ and $\delta
^{\prime}$
tend to zero. This implies that $V_{t,x}^{\epsilon,\delta}$ converges in
$L^{2}$, as $\varepsilon$ and $\delta$ tend to zero, to a limit
denoted by $%
V_{t,x}$. On the other hand, from the above computations, we have
\begin{eqnarray*}
E^W [ ({V_{t,x}^{\epsilon,\delta}} )^2 ]
&=&\alpha_{H_0} \int_{[0,t]^{4}} p_{2\epsilon}(B_{s}-B_{r})\varphi_{\delta}(t-s-u) \\
&&\hspace*{41.5pt}{}\times\varphi_{\delta}(t-r-v)
|u-v|^{2H_{0}-2}\,du\,dv\,ds\,dr
\end{eqnarray*}
and this expression
converges to right-hand side of
(\ref{e.8.9}) almost surely. Moreover, because of the above
arguments, the
convergence is also in $L^1$ and this implies (\ref{e.8.9}).
\end{pf}

%s6.3 ###
\subsection{Feynman--Kac formula}

By Proposition 3.3 and Theorem 6.2 in \cite{hunualart}, we have the
following theorem.
\begin{theorem}
Suppose that $H_{1}=1/2$ and $H_{0}>3/4$. Then, for any $\lambda\in
\mathbb{%
R}$, we have
\[
E\exp\biggl( \lambda\int_{0}^{t}\int_{\mathbb{R}}\delta
(B_{t-r}^{x}-y)W(dr,dy) \biggr) <\infty
\]
and, for any measurable and bounded function $f$, the process%
%
%e6.3 ###
%
\begin{equation} \label{utxa}
u(t,x)=E^{B} \biggl( f(B_{t}^{x})\exp\biggl( \int_{0}^{t}\int_{\mathbb{R}%
}\delta(B_{t-r}^{x}-y)W(dr,dy) \biggr) \biggr)
\end{equation}
is a weak solution of (\ref{e.1.1}).
\end{theorem}

%s6.4 ###
\subsection{\texorpdfstring{H\"{o}lder continuity}{Holder continuity}}

We also have the following theorem, whose proof is similar to that of
Theorem \ref{theo61}.
\begin{theorem}
Suppose that $H_1=1/2$, $H_0>3/4$ and let $u(t,x)$ be the solution of
(\ref{e.1.1}). Then, $u(t,x)$ has a continuous modification such that for
any $%
\rho\in(0,H_0-3/4) $ and any compact rectangle $I \subset\mathbb{R}%
_+\times\mathbb{R}$, there exists a positive random variable $K_I $
such that
almost surely, for any $(t_1,x_1),(t_2,x_2)\in I$, we have
\[
| u(t_2,x_2 )- u(t_1,x_1 )|\le K_I (|t_2-t_1|^\rho+|x_2-x_1|^{2\rho}).
\]
\end{theorem}
\begin{pf}
As in the proof of Theorem \ref{theo61}, we have
\[
E^W|u(s,x)-u(t,y)|^p\le C_p [E^B E^W|\tilde V-V|^2 ]^{p/2} ,
\]
where $V=\int_0^t \int_\RR{\delta}(B_{t-r}^x-z)W(dr,dz)$ and
$\tilde
V=\int_0^s \int_\RR{\delta}(B_{s-r}^y-z)W(dr,dz)$. If $s=t$, then we
can write
\begin{eqnarray*}
E^B E^W|\tilde V-V|^2 &=& 2 \int_0^t\int_0^t |r-v|^{2H_0-2} \\
&&\hspace*{35.3pt}{} \times E [\delta(B_{r}-B_{v})-\delta(B_{r }-B_{v}+x-y) ] \,dr\,dv \\
&=& \frac{2}{\sqrt{2\pi}}\int_0^t\int_0^t
|r-v|^{2H_0-5/2}\bigl(1-e^{-{
(x-y)^2}/({2|r-v|}) }\bigr)\,dr \,dv .
\end{eqnarray*}
For any $2\rho<{\gamma}<2H_0-3/2$, we have $ 1-e^{-{(x-y)^2%
}/({2 |r-v|}) }\le( \frac{(x-y)^{2 }}{2|r -v| } )^\ga$. Thus, $E^B
E^W|\tilde V-V|^2\le C_\ga|x-y|^{2{\gamma}}$. Consequently, we have
%
%e6.4 ###
%
\begin{equation} \label{e.7.10}
E^W|u(t,x)-u(t,y)|^p\le C |x-y|^{{\gamma} p}.
\end{equation}
On the other hand, if $x=y$,
then
\begin{eqnarray*}
&&E^B E^W|\tilde V-V|^2 \\
&&\qquad= C \biggl[ \int_{s} ^{t } \int_{s} ^{t}
|r-v|^{2H_0-5/2}\,dr\,dv\\
&&\qquad\quad\hspace*{11.7pt}{}+\int_0^{s}\int_0^{t}|r-v |^{2H_0-2} (|r-v|^{-1/2}-|r-v+
t-s|^{-1/2 } )\,dr\,ds \biggr]
\end{eqnarray*}
and, by a similar computation to step 3 before, we can obtain
%
%e6.5 ###
%
\begin{equation} \label{e.7.11}
E^W|u(s,x)-u(t ,x)|^p\le C (t-s)^{(H_0-3/4) p}.
\end{equation}
Combining (\ref{e.7.10}) and (\ref{e.7.11}), we prove the
theorem.
\end{pf}

%s6.5 ###
\subsection{Regularity of the density}

We can also show the following result.
\begin{theorem}
Suppose that $d=1$, $H_1=1/2$ and $H_0>3/4$. Fix $t>0$ and $x\in
\mathbb
{R} $. Assume that for any positive number $p$,
$E|f(B_{t}+x)|^{-p}<\infty$.
The law of $u(t,x)$ then has a smooth density.
\end{theorem}
\begin{pf}
The proof is similar to that of Theorem \ref{th2}, using the existence
of finite moments of all orders for the self-intersection local time of
the Brownian motion proved in the \hyperref[app]{Appendix} (see Proposition \ref{propA1}).
\end{pf}

%s7 ###
\section{Skorokhod-type equations and chaos expansion}\label{sec7}

In this section, we consider the following heat equation on $\mathbb{R}^d$:
%
%e7.1 ###
%
\begin{equation}\label{e.9.1}
\cases{\dfrac{\partial u}{\partial t}= \dfrac12 \Delta u+ u \diamond
\dfrac{\partial
^{d+1}}{\partial t\,\partial x_1 \cdots\partial x_d}W, \vspace*{2pt}\cr
u(0,x)=f(x).}
\end{equation}
The difference between the above equation and (\ref{e.1.1})
is that here we use the Wick product $\diamond$ (see, e.g., \cite
{hy}). This equation is studied in \cite{hunualart} in the case $%
H_1=\cdots=H_d=1/2$. As in that paper, we can define the following
notion of mild
solution.
\begin{definition}
\label{def1} An adapted random field $u=\{u(t,x),t\geq0,x\in\mathbb{R}
^{d}\}$ such that $E(u ^{2}(t,x))<\infty$ for all $(t,x)$ is a mild
solution to equation (\ref{e.9.1}) if, for any $(t,x)\in[
0,\infty
)\times\mathbb{R}^{d}$, the process $\{p_{t-s}(x-y)u(s,y)\mathbf{1}%
_{[0,t]}(s),s\geq0,y\in\mathbb{R}^{d}\mathbb{\}}$ is Skorokhod integrable
and the following equation holds:
%
%e7.2 ###
%
\begin{equation} \label{e.9.2}
u(t,x)=p_{t}f(x)+\int_{0}^{t}\int_{\mathbb
{R}^{d}}p_{t-s}(x-y)u(s,y)\,\delta
W_{s,y},
\end{equation}
where $p_t(x)$ denotes the heat kernel and $p_tf(x)=\int_{\mathbb{R}^d}
p_t(x-y) f(y)\,dy$.
\end{definition}

As in the paper \cite{hunualart}, the mild solution $u(t,x)$ of %
(\ref{e.9.1}) admits the following Wiener chaos expansion:
%
%e7.3 ###
%
\begin{equation} \label{e.9.3}
u(t,x)=\sum_{n=0}^{\infty}I_{n}(f_{n}(\cdot,t,x)),
\end{equation}
where $I_n$ denotes the multiple stochastic integral with respect to
$W$ and
$f_{n}(\cdot,t,x)$ is a symmetric element in $\mathcal{H} ^{\otimes n}$,
defined explicitly as
%
%e7.4 ###
%
\begin{eqnarray} \label{e.9.4}\qquad
&&f_{n}(s_{1},y_{1},\ldots,s_{n},y_{n},t,x)\nonumber\\[-8pt]\\[-8pt]
&&\qquad=\frac{1}{n!}
p_{t-s_{\sigma(n)}}\bigl(x-y_{\sigma(n)}\bigr)\cdots p_{s_{\sigma
(2)}-s_{\sigma(1)}}\bigl(y_{\sigma(2)}-y_{\sigma(1)}\bigr)p_{s_{\sigma
(1)}}f\bigl(y_{\sigma(1)}\bigr). \nonumber
\end{eqnarray}
In the above equation, $\sigma$ denotes a permutation of $\{1,2,\ldots
,n\}$
such that $0<s_{\sigma(1)}<\cdots<s_{\sigma(n)}<t$. Moreover, the
solution, if it exists, will be unique because
the kernels in the Wiener chaos expansion are uniquely determined.

The following theorem is the main result of this section.
\begin{theorem}
Suppose that $2H_0+\sum_{i=1}^d H_{i}>d+1$ and that $f$ is a bounded
measurable function. Then, the process%
%
%e7.5 ###
%
\begin{eqnarray}\label{e.9.5}\qquad
u(t,x) &=&E^{B} \Biggl[ f(B_{t}^{x})\exp\Biggl( \int_{0}^{t}\int_{\mathbb{R}%
^d}\delta(B_{t-r}^{x}-y)W(dr,dy) \nonumber\\
&&\hspace*{73.3pt}{} -\frac12 \alpha_H \int_{0}^{t}\int_{0}^{t}|r -s|^{2H_{0}-2}
\\
&&\qquad\quad\hspace*{105.2pt}{}\times\prod_{i=1}^d | B^i_{r }-B^i_{s} | ^{2H_{i}-2}\,dr \,ds \Biggr)
\Biggr]\nonumber
\end{eqnarray}
is the unique mild solution to equation (\ref{e.1.1}).
\end{theorem}
\begin{pf}
From Theorem \ref{teo2}, we obtain that the expectation $E^B$ in
(\ref
{e.9.5}) is well defined. It then suffices to show that the random variable
$u(t,x)$ has the Wiener chaos expansion (\ref{e.9.3}). This can be easily
proven by expanding the exponential and then taking the expectation with
respect to $B$.

Theorem \ref{teo1} implies that almost surely $\delta(B_{t-\cdot}^{x}-\cdot)$
is an element of $\mathcal{H}$ with a norm given by (\ref{e.3.4}).
As a consequence, almost surely with respect to the Brownian motion
$B$, we
have the following chaos expansion for the exponential factor in
equation (%
\ref{e.9.5}):
\begin{eqnarray*}
&& \exp\Biggl( \int_{0}^{t}\int_{\mathbb{R}^d}\delta(B_{t-r}^{x}-y)W(dr,dy)
\\
&&\qquad{}-\frac12 \alpha_H \int_{0}^{t}\int_{0}^{t}|r -s|^{2H_{0}-2} \prod_{i=1}^d
| B^i_{r }-B^i_{s} | ^{2H_{i}-2}\,dr \,ds \Biggr)
=\sum_{n=0}^\infty
I_n(g_n) ,
\end{eqnarray*}
where $g_n$ is the symmetric element in $\mathcal{H}^{\otimes n}$
given by
%
%e7.6 ###
%
\begin{equation}
g_n (s_{1},y_{1},\ldots,s_{n},y_{n},t,x)=\frac{1}{n!}{\delta}%
(B_{t-s_1}^x-y_1) \cdots{\delta}(B_{t-s_n}^x-y_n) .
\end{equation}
Thus, the right-hand side of (\ref{e.9.5}) admits the chaos
expansion
%
%e7.7 ###
%
\begin{equation}
u(t,x)=\sum_{n=0}^\infty\frac1{n!} I_n(h_n(\cdot,t,x))
\end{equation}
with
%
%e7.8 ###
%
\begin{equation}
h_n(t,x)=E^B [f(B_t^x) {\delta}(B_{t-s_1}^x-y_1) \cdots{\delta}%
(B_{t-s_n}^x-y_n) ] .
\end{equation}
This can be regarded as a Feynman--Kac formula for the coefficients of
the chaos
expansion of the solution of (\ref{e.9.1}). To compute the above
expectation, we shall use the following identity:
%
%e7.9 ###
%
\begin{eqnarray} \label{e.9.7}
E^B [f(B_t^x) {\delta}(B_t^x-y) | \mathcal{F}_s ] &=& \int_{%
\mathbb{R}^d} p_{t-s}(B_s^x -z) f(z) {\delta}(z-y) \,dz
\nonumber\\[-8pt]\\[-8pt]
&=& p_{t-s}(B_s^x -y)f(y) .\nonumber
\end{eqnarray}
Assume that $0< s_{{\sigma}(1)}<\cdots<s_{{\sigma}(n)}<t$ for some
permutation ${\sigma}$ of $\{1, 2, \ldots, n\}$. Then, conditioning with
respect to $\mathcal{F}_{t-s_{{\sigma}(1)}}$ and using the Markov property
of the Brownian motion, we have
\begin{eqnarray*}
h_n(t,x) & =& E^B \bigl\{ E^B \bigl[ {\delta}\bigl(B_{t-s_{{\sigma
}(n)}}^x-y_{{\sigma%
}(n)}\bigr) \cdots\\
&&\hspace*{17pt}{} \times{\delta}\bigl(B_{t-s_{{\sigma}(1)}}^x-y_{{\sigma}(1)}\bigr)
f(B_t^x) |\mathcal{F}_{t-s_{{\sigma}(1)}} \bigr] \bigr\} \\
& =& E^B \bigl[ {\delta}\bigl(B^x_{t-s_{{\sigma}(n)}} -y_{{\sigma}(n)} \bigr)
\cdots{%
\delta}\bigl(B_{t-s_{{\sigma}(1)} }^x-y_{{\sigma}(1)} \bigr) p_{s_{{\sigma}(1)}}
f\bigl(B_{t-s_{{\sigma}(1)}}^x \bigr) \bigr] .
\end{eqnarray*}
Conditioning with respect to $\mathcal{F}_{t-s_{{\sigma}(2)}}$ and
using (\ref{e.9.7}), we have
\begin{eqnarray*}
h_n(t,x) &=&E^B \bigl\{ E^B \bigl[ {\delta}\bigl(B_{t-s_{{\sigma}(n)}^x}-y_{{\sigma}
(n)} \bigr) \\
&&\hspace*{16.1pt}{}\times{\delta}\bigl(B_{t-s_{{\sigma}(1)} }^x-y_{{\sigma}(1)} \bigr)
p_{s_{{\sigma}%
(1)}} f\bigl(B_{t-s_{{\sigma}(1)}}^x \bigr) \bigr] |\mathcal{F}_{t-s_{{\sigma}(2)}}
\bigr\} \\
&=&E^B \bigl\{ {\delta}\bigl(B_{t-s_{{\sigma}(n)}^x}-y_{{\sigma}(n)} \bigr) \cdots{
\delta}\bigl(B_{t-s_{{\sigma}(2)} }^x -y_{{\sigma}(2)} \bigr) \\
&&\hspace*{16.7pt}{}\times E^B \bigl[ {\delta}\bigl(B_{t-s_{{\sigma}(1)} }^x-y_{{\sigma}(1)} \bigr)
p_{s_{{\sigma}(1)}} f\bigl(B_{t-s_{{\sigma}(1)}}^x \bigr) |\mathcal
{F}_{t-s_{{\sigma%
}(2)}} \bigr] \bigr\} \\
&=&E^B \bigl[ {\delta}\bigl(B_{t-s_{{\sigma}(n)}^x}-y_{{\sigma}(n)} \bigr) \cdots{%
\delta}\bigl(B_{t-s_{{\sigma}(2)} }^x-y_{{\sigma}(2)} \bigr) \\
&&\hspace*{15.9pt}{}\times p_{s_{{\sigma}(2)}- s_{{\sigma}(1)}}\bigl(B_{t-s_{{\sigma}(2)}
}^x-y_{{\sigma}(1)} \bigr) p_{ s_{{\sigma}(1)}} f\bigl(y_{{\sigma}(1)} \bigr) \bigr] .
\end{eqnarray*}
Continuing in this way, we find that
\[
h_n(t,x) =p_{t-s_{\sigma(n)}}\bigl(x-y_{\sigma(n)}\bigr)\cdots p_{s_{\sigma
(2)}-s_{\sigma(1)}}\bigl(y_{\sigma(2)}-y_{\sigma(1)}\bigr)p_{s_{\sigma
(1)}}f\bigl(y_{\sigma(1)}\bigr),
\]
which is the same as (\ref{e.9.4}).
\end{pf}
\begin{remark}
The method of this section can be applied to obtain a Feynman--Kac formula
for the coefficients of the chaos expansion of the solution of equation
(\ref%
{e.1.1}):
\[
u(t,x)=\sum_{n=0}^{\infty}\frac{1}{n!}I_{n}(h_{n}(\cdot,t,x))
\]
with
%
%e7.10 ###
%
\begin{eqnarray}\qquad
h_{n}(t,x) &=&E^{B} \Biggl[f(B_{t}^{x}){\delta}(B_{t-s_{1}}^{x}-y_{1})\cdots
{%
\delta}(B_{t-s_{n}}^{x}-y_{n}) \nonumber\\
&&\hspace*{18.7pt}{} \times\exp\Biggl( \frac{1}{2}\alpha
_{H}\int_{0}^{t}\int_{0}^{t}|r-s|^{2H_{0}-2}\\
&&\hspace*{106.2pt}{}\times\prod_{i=1}^{d} |
B_{r}^{i}-B_{s}^{i} | ^{2H_{i}-2}\,dr\,ds \Biggr) \Biggr] .\nonumber
\end{eqnarray}
\end{remark}
\begin{remark}
We can also consider equation (\ref{e.1.1}) when $d=1$, $H_{1}=1/2$ and
$%
H_{0}>3/4$. In this case, we easily see that the solution $u(t,x)$
admits the
following chaos expansion:
\[
u(t,x)=\sum_{n=0}^{\infty}\frac{1}{n!}I_{n}(h_{n}(\cdot,t,x))
\]
with
%
%e7.11 ###
%
\begin{eqnarray}\qquad
h_{n}(t,x) &=&E^{B} \biggl[f(B_{t}^{x}){\delta}(B_{t-s_{1}}^{x}-y_{1})\cdots
{%
\delta}(B_{t-s_{n}}^{x}-y_{n}) \nonumber\\
&&\hspace*{18.1pt}{} \times\exp\biggl( \frac{1}{2}\alpha
_{H_{0}}\int_{0}^{t}\int_{0}^{t}|r-s|^{2H_{0}-2}\delta(
B_{r}-B_{s} ) \,dr\,ds \biggr) \biggr] .
\end{eqnarray}
\end{remark}

From the Feynman--Kac formula, we can derive the following formula for
the moments of the solution analogous to (\ref{e.4.7}), which can be
compared with the formulas obtained in \cite{hunualart} in the case
$H_1= \cdots=H_d=\frac12$:
\begin{eqnarray*}
&&E ( u(t,x)^{p} )\\
&&\qquad=E \Biggl(\prod_{j=1}^{p}f(B_{t}^{j}+x)
\\
&&\qquad\quad\hspace*{30pt}{}\times\exp\Biggl[ \alpha_{H} \sum_{j,k=1, j<k }^{p}\int_{0}^{t}%
\int_{0}^{t}|s-r|^{2H_{0}-2}\\
&&\qquad\quad\hspace*{150pt}{}\times\prod_{i=1}^{d}
|B_{s}^{j,i}-B_{r}^{k,i}|^{2H_{i}-2}\,ds\,dr \Biggr] \Biggr),
\end{eqnarray*}
where $p\ge1$ is an integer and $B^j, 1\le j\le d$, are independent
$d$-dimensional Brownian motions.

\begin{appendix}\label{app}
%s8 ###
\section*{Appendix}

\setcounter{lemma}{0}
\begin{lemma}
\label{lemma1} Suppose that $0<\alpha<1$, $\epsilon>0$, $x>0$ and that
$X $ is a
standard normal random variable. Then, there is a constant $C$,
independent of
$x$ and $\epsilon$ (it may depend on ${\alpha}$), such that
\[
E|x+\epsilon X|^{-\alpha}\le C\min(\epsilon^{-\alpha},x^{-\alpha}) .
\]
\end{lemma}
\begin{pf}
It is straightforward to check that $K=\sup_{z\ge0}E|z+X|^{-\alpha
}<\infty$.
Thus,
%
%e8.1 ###
%
\setcounter{equation}{0}
\begin{equation} \label{e.2.1}
E|x+\epsilon X|^{-\alpha}=\epsilon^{-\alpha}E\biggl|\frac{x}{\epsilon}%
+X\biggr|^{-\alpha}\le K\epsilon^{-\alpha} .
\end{equation}
On the other hand,
\begin{eqnarray*}
E|x+\epsilon X|^{-\alpha} &=&\frac{1}{\sqrt{2\pi}} \int_\mathbb
{R}%
|x+\epsilon y|^{-\alpha}e^{-{y^2}/{2}}\,dy \\
&=&\frac{1}{\sqrt{2\pi}} \biggl(\int_{\{|x+\epsilon y|> {x}/{2}\}}
|x+\epsilon y|^{-\alpha} e^{-{y^2}/{2}}\,dy \\
&&\hspace*{30pt}{} +\int_{\{|x+\epsilon y|\le{x}/{2}\}} |x+\epsilon y|^{-\alpha
}e^{-%
{y^2}/{2}}\,dy \biggr) .
\end{eqnarray*}
It is easy to see that the first integral is bounded by $Cx^{-{\alpha
}}$ for
some constant $C$. The second integral, denoted by $B$, is bounded as
follows:
\begin{eqnarray*}
B&=&C\frac{1}{\epsilon}\int_{|z|<{x/2}}|z|^{-\alpha
}e^{-
{(z-x)^2%
}/({2\epsilon^2})}\,dz \le C\frac{1}{\epsilon}\int_{|z|<{x/2}%
}|z|^{-\alpha}e^{-{x^2}/({8\epsilon^2})}\,dz \\
&=&C\frac{x}{\epsilon}e^{-{x^2}/({8\epsilon^2})}x^{-\alpha} \le
Cx^{-\alpha} .
\end{eqnarray*}
Thus, we have $E|x+\epsilon X|^{-\alpha}\le C |x|^{-{\alpha}}$. Combining
this with (\ref{e.2.1}), we obtain the lemma.
\end{pf}
\begin{lemma}
\label{lemma2} Suppose that $\alpha\in(0,1)$. There exists a
constant $C>0$
such that
\[
\sup_{\epsilon,\epsilon^{\prime}} \int_{\mathbb{R}^2}p_\epsilon(x_1+y_1)
p_{\epsilon^{\prime}}(x_2+y_2)|y_1-y_2|^{-\alpha}\,dy_1\,dy_2\le
C|x_1-x_2|^{-\alpha}.
\]
\end{lemma}
\begin{pf}
We can write
\[
\int_{\mathbb{R}^2}p_\epsilon(x_1+y_1)p_{\epsilon^{%
\prime}}(x_2+y_2)|y_1-y_2|^{-\alpha}\,dy_1\,dy_2= { E} ( |{\varepsilon}
X_1-x_1-{\varepsilon}^{\prime}X_2+x_2|^{-{\alpha}} ) .
\]
Thus, Lemma \ref{lemma2} follows directly from Lemma \ref{lemma1}.
\end{pf}
\begin{lemma}
\label{lemma3} Suppose that $\alpha\in(0,1)$. There exists a
constant $C>0$
such that
\[
\sup_{\delta,\delta^{\prime}}\int_0^t\int_0^t\varphi_\delta
(t-s_1-r_1)%
\varphi_{\delta^{\prime}}(t-s_2-r_2)|r_1-r_2|^{-\alpha}\,dr_1\,dr_2 \le
C|s_1-s_2|^{-\alpha}.
\]
\end{lemma}
\begin{pf}
Since
\[
p_\delta(x)\ge p_\delta(x) I_{[0,\sqrt\delta]}(x)=\frac{1}{\sqrt
{2\pi
\delta%
}}e^{-{x^2}/({2\delta})}I_{[0,\sqrt\delta]}(x)\ge\frac
{1}{\sqrt
{2\pi e }}%
\varphi_{\sqrt\delta}(x) ,
\]
the lemma follows from Lemma \ref{lemma2}.
\end{pf}
\begin{lemma}
\label{lemma4} Suppose that $2H_{0}+\sum_{i=1}^{d}H_{i}>d+1$. Let $%
B^{1},\ldots,B^{d}$ be independent one-dimensional Brownian motions. We then
have
\[
E \Biggl(
\int_{0}^{t}s^{2H_{0}-2}\prod_{i=1}^{d}|B_{s}^{i}|^{2H_{i}-2}\,ds \Biggr)
^{2}<\infty.
\]
\end{lemma}
\begin{pf}
We can write
\begin{eqnarray*}
&&E \Biggl(
\int_{0}^{t}s^{2H_{0}-2}\prod_{i=1}^{d}|B_{s}^{i}|^{2H_{i}-2}\,ds \Biggr)
^{2}\\
&&\qquad=2\int_{0}^{t}\int_{0}^{s}(sr)^{2H_{0}-2}
\prod_{i=1}^{d}E(|B_{s}^{i}|^{2H_{i}-2}|B_{r}^{i}|^{2H_{i}-2})\,dr\,ds.
\end{eqnarray*}
Let $X$ be a standard normal random variable. From Lemma \ref{lemma1},
taking into account that $2-2H_{i}<1$, we have, when $r<s$,
that
%
%e8.2 ###
%
\begin{eqnarray} \label{e.2.3}\qquad
E(|B_{r}^{i}|^{2H_{i}-2}|B_{s}^{i}|^{2H_{i}-2})
&=&E\bigl[\bigl|B_{r}^{i}|^{2H_{i}-2}E\bigl[|\sqrt{s-r}X+x\bigr|^{2H_{i}-2}|_{x=B_{r}^{i}}\bigr]\bigr]
\nonumber\\
&\leq&CE\bigl[|B_{r}^{i}|^{2H_{i}-2}(s-r)^{H_{i}-1}\bigr) \\
&\leq&Cr^{H_{i}-1}(s-r)^{H_{i}-1} .\nonumber
\end{eqnarray}
As a consequence, the conclusion of the lemma follows from the fact that
\[
\int_{0}^{t}\int_{0}^{s}r^{2H_{0}+%
\sum_{i=1}^{d}H_{i}-d-2}s^{2H_{0}-2}(s-r)^{\sum
_{i=1}^{d}H_{i}-d}\,dr\,ds<\infty,
\]
because $2H_{0}+\sum_{i=1}^{d}H_{i}-d-2>-1$ and $\sum_{i=1}^{d}H_{i}-d>-1$.
\end{pf}
\begin{lemma}
\label{lemma5} Let $B^{1},\ldots,B^{d}$ be independent one-dimensional
Brownian motions. If $\alpha_{i}\in(-1,0)$, $i=1,\ldots, d$, and $%
\sum_{i=1}^{d}\alpha_{i}>-2$, then
\[
E\exp\Biggl( \lambda
\int_{0}^{1}\prod_{i=1}^{d}|B_{s}^{i}|^{\alpha_{i}}\,ds \Biggr) <\infty
\]
for all $\lambda>0$.
\end{lemma}
\begin{pf}
The proof is based on the method of moments. We can write
\begin{eqnarray*}
E\exp\Biggl(\lambda\int_0^1\prod_{i=1}^d|B_s^i|^{\alpha_i}\,ds \Biggr) &=&
\sum_{n=1}^\infty\frac{\lambda^n}{n!} E\int_{[0,1]^n}
\prod_{k=1}^n\prod_{i=1}^d|B_{s_k}^i|^{\alpha_i}\,ds \\
&=&\sum_{n=1}^\infty\lambda^n \int_{[0<s_1<\cdots<s_n<1]} \prod_{i=1}^d
E \Biggl(\prod_{k=1}^n|B_{s_k}^i|^{\alpha_i} \Biggr)\,ds.
\end{eqnarray*}
From Lemma \ref{lemma1}, since $\alpha_i \in(-1,0)$, we obtain
\[
E [|B^i_{s_k}|^{\alpha_i}|\mathcal{F}^i_{s_{k-1}} ]=E [%
|B^i_{s_k}-B^i_{s_{k-1}}+B^i_{s_{k-1}}|^{\alpha_i}|\mathcal
{F}^i_{s_{k-1}}%
]\le C(s_k-s_{k-1})^{\alpha_i/2},
\]
where $\mathcal{F}_t$ is the filtration generated by the Brownian
motion $%
B^i $. As a consequence, taking the conditional expectation of $%
\prod_{k=1}^n|B_{s_k}^i|^{\alpha_i}$ with respect to the $\sigma
$-fields $%
\mathcal{F}^i_{s_{n-1}}, \mathcal{F}^i_{s_{n-2}}, \ldots, \mathcal{F}
^i_{s_{1}}$ and $\mathcal{F}^i_0$, we get
\[
E \Biggl(\prod_{k=1}^n|B_{s_k}^i|^{\alpha_i} \Biggr)\le C^n
(s_n-s_{n-1})^{\alpha_i/2}\cdots(s_2-s_{1})^{\alpha_i/2} s_1^{\alpha_i/2}.
\]
Letting $\alpha=\sum_{i=1}^d \alpha_i$, we have
\begin{eqnarray*}
&&E\exp\Biggl(\lambda\int_0^1\prod_{i=1}^d|B_s^i|^{\alpha_i}\,ds
\Biggr)\\
&&\qquad\le\sum_{n=1}^\infty(C\lambda)^n
\int_{[0<s_1<\cdots<s_n<1]}
(s_n-s_{n-1})^{\alpha/2}\cdots\\
&&\qquad\quad\hspace*{109pt}{}\times(s_2-s_{1})^{\alpha/2} s_1^{\alpha/2}\,ds.
\end{eqnarray*}
Since $\alpha> -2$, the integrals on the right-hand side are equal to
$\frac{%
(\Gamma(\alpha/2+1) )^n}{(n+n\alpha/2)\Gamma(n+n\alpha/2)}$ and
the series converges for any $\lambda>0$.
\end{pf}
\begin{lemma}
\label{Lem7.6} For any $0<\alpha<1$, define
\[
C_ \alpha(y)=E (|\xi|^{-\alpha}-|y+\xi|^{-\alpha} ),
\]
where $y>0$ and $\xi$ is a standard normal random variable. Then,
\[
C_\alpha(y) \le C \min\bigl(1, ( y^2 + y^{3-\alpha}) \bigr)
\]
for some constant $C>0$.
\end{lemma}
\begin{pf}
First, note that $C_\alpha(y)<C$, where $C>0$ is a constant, since $%
\lim_{y\to\infty}E|y+\xi|^{-\alpha}=0$. On the other hand, we can decompose
the function $C_\alpha(y)$ as follows:
\begin{eqnarray*}
C_\alpha(y) &=&\frac{1}{\sqrt{2\pi}}\int_{\mathbb{R} } (%
|x|^{-\alpha}-|y+x|^{-\alpha} )e^{-{x^2}/{2}}\,dx \\
&=&\frac{1}{\sqrt{2\pi}} \biggl(\int_{\{x\ge0\} \cup\{x\le-y\}} (%
|x|^{-\alpha}-|y+x|^{-\alpha} )e^{-{x^2}/{2}}\,dx \\
&&\hspace*{35.9pt}{} +\int_{\{-y<x<0\}} (|x|^{-\alpha}-|y+x|^{-\alpha} )e^{-{%
x^2}/{2}}\,dx \biggr) \\
&=&\frac{1}{\sqrt{2\pi}}(A+B) ,
\end{eqnarray*}
where $A$ and $B$ denote the first and second integrals, respectively, in
the second-to-last line. For integral $A$, we can write
\begin{eqnarray*}
A&=&\int_0^\infty\bigl(x^{-\alpha}-(x+y)^{-\alpha} \bigr) \bigl(e^{-x^2/2}-
e^{-(x+y)^2/2} \bigr)\,dx \\
&\le&\int_0^\infty x^{-\alpha}(x+y)^{1-\alpha} [(x+y)^\alpha-
x^\alpha] y
e^{-{x^2}/{2}}\,dx.
\end{eqnarray*}
Therefore,
\begin{eqnarray*}
A &\le& \int_0^\infty x^{1-2\alpha} [(x+y)^\alpha- x^\alpha] y
e^{-{x^2}/{2}}\,dx\\
&&{} +\int_0^\infty x^{-\alpha} [(x+y)^\alpha- x^\alpha]
y^{2-\alpha}
e^{-{x^2}/{2}}\,dx.
\end{eqnarray*}
For the first integral in the above expression, we use the estimate $%
(x+y)^\alpha- x^\alpha\le\alpha y x^{ \alpha-1}$ and for the
second, we
use $(x+y)^\alpha- x^\alpha\le y^{ \alpha}$. In this way, we obtain
\[
A\le Cy^2
\]
for some constant $C>0$. On the other hand,
\[
B = \int_0^y x^{-\alpha}\bigl(e^{-{x^2}/{2}}-e^{-{(x+y)^2}/{2}}\bigr)\,dx
\le
\int_0^y x^{-\alpha}(x+y)y\,dx \le Cy^{3-\alpha}
\]
for some constant $C>0$, which completes the proof of the lemma.
\end{pf}
\begin{prop} \label{propA1}
Let $B$ be a one-dimensional standard Brownian motion. Then, for any $p>0$,
\[
E \biggl|\int_0^1\int_0^1 \delta(B_t-B_s)\,ds\,dt \biggr|^{-p}<\infty.
\]
\end{prop}
\begin{pf}
For $k=1,\ldots, 2^{n-1}$, we define $%
A_{n,k}= [ \frac{2k-2}{2^{n}},\frac{2k-1}{2^{n}} ]
\times[ \frac{2k-1}{2^{n}},\frac{2k}{2^{n}} ] $ and
\[
\alpha_{n,k}=\int_{A_{n,k}}\delta(B_t- B_s)\,ds\,dt.
\]
The random variables $\alpha_{n,k}$ have the following two
properties:
\begin{longlist}
\item for every $n\ge1$, the variables $\alpha_{n,1},\ldots,%
\alpha_{n,2^{n-1}}$ are independent;

\item $\alpha_{n,k}\stackrel{d}{=}2^{-n/2}
\int_0^1\int_0^1\delta(B_t-\tilde B_s)\,ds\,dt$ and $\tilde B $ is a
standard Brownian motion independent of $B$.
\end{longlist}
For any $p>0$, we may choose a integer $n>0$ such that
$p2^{1-n}<1/3$. Then, we can write
\[
E \biggl|\int_0^1\int_0^1
\delta(B_t-B_s)\,ds\,dt \biggr|^{-p} \le
E \Biggl|\sum_{k=1}^{2^{n-1}}\alpha_{n,k} \Biggr|^{-p}\le
E \Biggl|\prod_{k=1}^{2^{n-1}}\alpha_{n,k} \Biggr|^{-p2^{1-n}}
\]
and it suffices to show that $E |\int_0^1\int_0^1\delta(B_t-\tilde
B_s)\,ds\,dt |^{-p}<\infty$
for some $p>0$. Notice that
\[
L:=\int_0^1\int_0^1\delta(B_t-\tilde
B_s)\,ds\,dt=\int_ {\mathbb{R}} L_1^x\tilde
L_1^x\,dx,
\]
where $L_t^x$ (resp., $L_t^ {\tilde{x}}$) denotes the local time of the
Brownian motion $B$ (resp., $\tilde B$).
As a consequence, for any $0<\alpha<1$,
\begin{eqnarray*}
P(L<\epsilon) & \le& P \biggl(\int_0^{ \epsilon^{4/5}} L_1^x\tilde
L_1^x\,dx \biggr) \\
&\le& P \biggl( L_1^0\tilde L_1^0 <\frac12 \epsilon^{1/5} \biggr)
+P \biggl( \int_0^{ \epsilon^{4/5}} | L_1^0\tilde L_1^0-L_1^x\tilde L_1^x|
\,dx \ge\frac\epsilon2 \biggr) \\
&\le& \frac1 {\sqrt{2}} \epsilon^{1/10} (E( L^0_1) ^{-1/2})^2 +
\frac
2\epsilon
\int_0^{ \epsilon^{4/5}}E | L_1^0\tilde L_1^0-L_1^x\tilde L_1^x|\, dx \\
&\le& C \epsilon^{1/10},
\end{eqnarray*}
which implies that $E(L^{1/10}) <\infty$.
\end{pf}
\end{appendix}

\section*{Acknowledgments}
We would like to thank an anonymous referee for his/her valuable
comments which helped us to improve this paper.

% imsref loaded by lrinkeviciute, 2010-08-03 07:41:10
%

%
\printaddresses


\begin{thebibliography}{19}

%%b1 ###
%%
%(\byear{1999}).
%In \bbooktitle{Stochastic Analysis, Control, Optimization and Applications}.
%%\bseries{Systems Control Found. Appl.}
%%
%
%%b2 ###
%%
%(\byear{1995}).
%stochastic partial differential equations}.
%%
%
%%b3 ###
%%
%(\byear{1998}).
%model with
%continuous space parameter}.
%%
%
%%b4 ###
%%
%(\byear{2007}).
%%\bseries{Chapman \& Hall/CRC Applied Mathematics and Nonlinear Science Series}.
%%

%b5 ###
\bibitem{DS}
%
\begin{barticle}[mr]
\bauthor{\bsnm{Dawson},~\bfnm{Donald~A.}\binits{D.~A.}} \AND
\bauthor{\bsnm{Salehi},~\bfnm{Habib}\binits{H.}}
(\byear{1980}).
\btitle{Spatially homogeneous random evolutions}.
\bjournal{J. Multivariate Anal.}
\bvolume{10}
\bpages{141--180}.
\bid{doi={10.1016/0047-259X(80)90012-3}, mr={575923}}
\end{barticle}
%
\endbibitem

%b6 ###
\bibitem{fred}
%
\begin{bbook}[mr]
\bauthor{\bsnm{Freidlin},~\bfnm{Mark}\binits{M.}}
(\byear{1985}).
\btitle{Functional Integration and Partial Differential Equations}.
\bseries{Annals of Mathematics Studies}
\bvolume{109}.
\bpublisher{Princeton Univ. Press}, \baddress{Princeton, NJ}.
\bid{mr={833742}}
\end{bbook}
%
\endbibitem

%b7 ###
\bibitem{Hi}
%
\begin{bmisc}[vtex]
\bauthor{\bsnm{Hinz},~\bfnm{H.}\binits{H.}}
(\byear{2009}).
\bhowpublished{Burgers system with a fractional Brownian
random force. Preprint,
Technische Univ. Berlin}.
\end{bmisc}
%
\endbibitem

%b8 ###
\bibitem{HLN}
%
\begin{bmisc}[vtex]
\bauthor{\bsnm{Hu},~\bfnm{Y.}\binits{Y.}},
\bauthor{\bsnm{Lu},~\bfnm{F.}\binits{F.}} \AND
\bauthor{\bsnm{Nualart},~\bfnm{D.}\binits{D.}}
(\byear{2010}).
\bhowpublished{Feynman--Kac formula for the heat equation
driven by fractional noise with Hurst parameter $H<1/2$.
Preprint, Univ. Kansas}.
\end{bmisc}
%
\endbibitem

%b9 ###
\bibitem{hunualart}
%
\begin{barticle}[mr]
\bauthor{\bsnm{Hu},~\bfnm{Yaozhong}\binits{Y.}} \AND
\bauthor{\bsnm{Nualart},~\bfnm{David}\binits{D.}}
(\byear{2009}).
\btitle{Stochastic heat equation driven by fractional noise and local time}.
\bjournal{Probab. Theory Related Fields}
\bvolume{143}
\bpages{285--328}.
\bid{doi={10.1007/s00440-007-0127-5}, mr={2449130}}
\end{barticle}
%
\endbibitem

%b10 ###
\bibitem{hy}
%
\begin{barticle}[vtex]
\bauthor{\bsnm{Hu},~\bfnm{Yao-Zhong}\binits{Y.-Z.}} \AND
\bauthor{\bsnm{Yan},~\bfnm{Jia-An}\binits{J.-A.}}
(\byear{2009}).
\btitle{Wick calculus for nonlinear {G}aussian functionals}.
\bjournal{Acta Math. Appl. Sin. Engl. Ser.}
\bvolume{25}
\bpages{399--414}.
\bid{doi={10.1007/s10255-008-8808-0}, mr={2506982}}
\end{barticle}
%
\endbibitem

%b11 ###
\bibitem{kunita}
%
\begin{bbook}[mr]
\bauthor{\bsnm{Kunita},~\bfnm{Hiroshi}\binits{H.}}
(\byear{1990}).
\btitle{Stochastic Flows and Stochastic Differential Equations}.
\bseries{Cambridge Studies in Advanced Mathematics}
\bvolume{24}.
\bpublisher{Cambridge Univ. Press}, \baddress{Cambridge}.
\bid{mr={1070361}}
\end{bbook}
%
\endbibitem

%b12 ###
\bibitem{legall}
%
\begin{bincollection}[mr]
\bauthor{\bsnm{Le~Gall},~\bfnm{Jean-Fran{\c{c}}ois}\binits{J.-F.}}
(\byear{1994}).
\btitle{Exponential moments for the renormalized self-intersection
local time
of planar {B}rownian motion}.
In \bbooktitle{S\'eminaire de {P}robabilit\'es, {XXVIII}}.
\bseries{Lecture Notes in Math.}
\bvolume{1583}
\bpages{172--180}.
\bpublisher{Springer}, \baddress{Berlin}.
\bid{doi={10.1007/BFb0073845}, mr={1329112}}
\end{bincollection}
%
\endbibitem

%b13 ###
\bibitem{movi}
%
\begin{barticle}[mr]
\bauthor{\bsnm{Mocioalca},~\bfnm{Oana}\binits{O.}} \AND
\bauthor{\bsnm{Viens},~\bfnm{Frederi}\binits{F.}}
(\byear{2005}).
\btitle{Skorohod integration and stochastic calculus beyond the fractional
{B}rownian scale}.
\bjournal{J. Funct. Anal.}
\bvolume{222}
\bpages{385--434}.
\bid{doi={10.1016/j.jfa.2004.07.013}, mr={2132395}}
\end{barticle}
%
\endbibitem

%b14 ###
\bibitem{nualart}
%
\begin{bbook}[mr]
\bauthor{\bsnm{Nualart},~\bfnm{David}\binits{D.}}
(\byear{2006}).
\btitle{The {M}alliavin Calculus and Related Topics},
\bedition{2nd} ed.
\bpublisher{Springer}, \baddress{Berlin}.
\bid{mr={2200233}}
\end{bbook}
%
\endbibitem

%b15 ###
\bibitem{RV}
%
\begin{barticle}[mr]
\bauthor{\bsnm{Russo},~\bfnm{Francesco}\binits{F.}} \AND
\bauthor{\bsnm{Vallois},~\bfnm{Pierre}\binits{P.}}
(\byear{1993}).
\btitle{Forward, backward and symmetric stochastic integration}.
\bjournal{Probab. Theory Related Fields}
\bvolume{97}
\bpages{403--421}.
\bid{doi={10.1007/BF01195073}, mr={1245252}}
\end{barticle}
%
\endbibitem

%%b16 ###
%%
%(\byear{2002}).
%manifold}.
%%
%
%%b17 ###
%%
%(\byear{2007}).
%sub-{$n$}th chaos processes}.
%%

%b18 ###
\bibitem{VZ}
%
\begin{barticle}[mr]
\bauthor{\bsnm{Viens},~\bfnm{Frederi~G.}\binits{F.~G.}} \AND
\bauthor{\bsnm{Zhang},~\bfnm{Tao}\binits{T.}}
(\byear{2008}).
\btitle{Almost sure exponential behavior of a directed polymer in a fractional
{B}rownian environment}.
\bjournal{J. Funct. Anal.}
\bvolume{255}
\bpages{2810--2860}.
\bid{doi={10.1016/j.jfa.2008.06.020}, mr={2464192}}
\end{barticle}
%
\endbibitem

%b19 ###
\bibitem{Wa}
%
\begin{bincollection}[vtex]
\bauthor{\bsnm{Walsh},~\bfnm{John~B.}\binits{J.~B.}}
(\byear{1986}).
\btitle{An introduction to stochastic partial differential equations}.
In \bbooktitle{\'{E}cole D'\'et\'e de Probabilit\'es de {S}aint--{F}lour,
{XIV}---1984}.
\bseries{Lecture Notes in Math.}
\bvolume{1180}
\bpages{265--439}.
\bpublisher{Springer}, \baddress{Berlin}.
\bid{mr={876085}}
\end{bincollection}
%
\endbibitem

\end{thebibliography}
\end{document}